\documentclass[a4paper]{article}
\usepackage{amsthm,amssymb,amsmath,enumerate,graphics,graphicx,psfrag}

\newtheorem{definition}{Definition}[subsection]
\newtheorem{proposition}[definition]{Proposition}
\newtheorem{theorem}[definition]{Theorem}
\newtheorem{question}[definition]{Question}
\newtheorem{corollary}[definition]{Corollary}
\newtheorem{lemma}[definition]{Lemma}
\newtheorem{problem}[definition]{Problem}
\newtheorem{example}[definition]{Example}
\newtheorem{conjecture}[definition]{Conjecture}

\newcommand{\kreis}[1]{\mathaccent"7017\relax #1}

\newcommand{\comment}[1]{}
\newcommand{\N}{\mathbb N}

\newcommand{\Z}{\mathbb Z}
\newcommand{\df}{d_{e/v}}

\newcommand{\dv}{d_v}
\newcommand{\de}{d_e}

\title{Extremal Infinite Graph Theory}
\author{Maya Stein\footnote{Supported by Fondecyt grant no.~11090141.}}
\date{28.12.2010}

\begin{document}
\maketitle

\begin{abstract}
We survey various aspects of infinite extremal graph theory and prove several new results. 
The lead role play the parameters  connectivity and degree. This includes the {\em end degree}.
Many open problems are suggested.
\end{abstract}

\vskip1cm

\section{Introduction}
\subsection{A short overview}
Until now, extremal graph theory usually meant {\em finite} extremal graph theory. 
New notions, as the {\em end degrees}~\cite{degree,hcs}, {\em circles} and {\em arcs}, and the topological viewpoint~\cite{DBook}, make it possible to create the infinite counterpart of the theory.
 We attempt here to give an overview of results and open problems that fall into this emerging area of infinite graph theory.

The paper divides into three parts. The first part is about forcing substructures with assumptions on the degree. We use the vertex-/edge-degree of the ends (for a definition see below) to force highly connected subgraphs and grid minors. For ensuring large complete minors, the vertex-/edge-degree is not enough, and we introduce a new notion, the relative degree, which accomplishes the task, at least for locally finite graphs. Related problems will be addressed along the way.

The second part is on  minimal higher connectivity and edge-connectivity of graphs, that is $k$-(edge-)connectivity for some $k\in\N$. This includes minimality with respect to edge deletion,  with respect to vertex deletion, and with resepct to taking subgraphs. The main questions here are the existence of vertices or ends of small degree, and bounds on the number of these. We also discuss whether minimal $k$-connected subgraphs in the sense(s) above exist in every $k$-connected graph. This will lead to a discussion of the problem in certain subspaces of the topological space associated to the graph.

The third and last part of the survey is on circles and arcs. These are the topologically defined analogues of cycles and paths in infinite graphs (see Section~\ref{sec:def}). We first discuss results and problems related to Hamilton circles, then move on to (topological) tree-packing and arboricity, and finally close the paper with a discussion of problems related to
 connectivity-preserving arcs and circles.

\subsection{Structure and degree}
Extremal graph theory in its strictest sense is all about forcing some palpable properties of a graph, very often some interesting substructure, by making assumptions on the overall density of the graph, conveniently expressed in terms of global parameters such as the average or minimum degree. Tur\'an's well-known theorem is a classical result in this direction, for a discussion of its extension to infinite graphs see~\cite{BBModern,mayaMCW2010}.

Another typical and important result in finite extremal graph theory, which can be found in any standard textbook, is the following theorem of Kostochka. It says that large average degree forces a large complete minor (and the function $f_1(k)$ is essentially the best possible bound~\cite{thomason}). 
\setcounter{subsection}{1}
\begin{theorem}$\!\!${\bf\cite{DBook}}\label{thm:minor}
There is a constant $c_1$ so that, for every $k\in\N$, if $G$ is a finite graph of average degree at least $f_1(k):=c_1 k\sqrt{\log k}$, then $G$ has a complete minor of order $k$.
\end{theorem} 

Also large topological minors can be forced with similar assumptions in finite graphs, as the following result, due to Bollob\'as and Thomason, states. 

\begin{theorem}\label{thm:topminor}$\!\!${\bf \cite{DBook}}
There is a constant $c_2$ so that, for every $k\in\N$, if $G$ is a finite graph of average degree at least $f_2(k):=c_2 k^2$, then $G$ has a complete topological minor of order $k$.
\end{theorem} 

Let us see how these results could extend to infinite graphs.
First of all we have to note that  it is not clear what the average degree of an infinite graph should be. We shall thus stick to the minimal degree as our `density-indicating' parameter. A minor, on the other hand, is defined in same way as for finite graphs, only that the branch-sets may now be infinite.\footnote{As long as our minors are locally finite, however (which will always be the case in this paper), it does not make any difference whether we allow infinite branch-sets or not. It is easy to see that any infinite branch-set of a locally finite minor may be restricted to a finite one.}

In rayless graphs we will then get a verbatim extension of Theorems~\ref{thm:minor} and~\ref{thm:topminor} (namely Theorem~\ref{prop:rayless}). This theorem will follow from a useful reduction theorem (Theorem~\ref{thm:reduRayless}), which states that every rayless graph of minimum degree $k$ has a finite subgraph of minimum degree $k$. These results will be presented in Section~\ref{sec:rayless}. 

In graphs with rays, however,  large minimal degree at the vertices is too weak to force any interesting substructure. This is so because  infinite trees may have arbitrarily large degrees, but they do not even have any $2$-connected subgraphs. 
So at first sight, our goal seems unreachable. At second thought, however, the example of the infinite tree just shows that we did not translate the term `large local densities' in the right way to infinite graphs. Only having every finite part of an infinite graph send out a large number of edges will not produce large overall density, if we do not require something to `come back' from infinity.

The most natural way to do this is to impose a condition on the ends of the graph. Ends are defined as the equivalence classes of rays (one-way infinite paths), under the equivalence relation of not being separable by any finite set of vertices. Ends have a long history, see~\cite{endsBerniElmar}.

In~\cite{degree} and in~\cite{hcs}, {\em end degrees} were introduced. In fact, two notions have turned out useful (for different purposes): the vertex-degree and the edge-degree of an end $\omega$. The {\em vertex-degree} of $\omega$ is defined as the maximum cardinality of a set of (vertex-)disjoint rays in $\omega$, and the {\em edge-degree}  is defined as the maximum cardinality of a set of edge-disjoint rays in $\omega$. These maxima exist~\cite{halin65}. 

Do these notions help to force density in infinite graphs? To some extent they do: A large minimum degree at the vertices together with a large minimum vertex-/edge-degree at the ends implies a certain dense substructure, which takes the form of a highly connected or edge-connected subgraph. 

More precisely, there is a function $f_v$ such that every graph of minimum degree resp.~vertex-degree $f_v(k)$ at the vertices and the ends has a $k$-connected subgraph, and there is also  a function $f_e$ such that every graph of minimum degree/edge-degree $f_e(k)$ at the vertices and the ends has a $k$-edge-connected subgraph. While $f_e$ is linear, $f_v$ is quadratic, and this is almost best possible. All these results are from~\cite{hcs} and will be presented in Section~\ref{sec:hcs}.

Related results will be discussed in Sections~\ref{sec:grid} and~\ref{sec:Vtrans}. In the latter, we shall see that in locally finite vertex-transitive graphs, $k$-connectivity is implied by much weaker assumptions. In fact, the $k$-(edge-)connectivity of a locally finite vertex-transitive graph is equivalent to all its ends having vertex-(resp.~edge-) degree $k$.  In Section~\ref{sec:grid} we shall see that independently of the degrees at the vertices, large vertex-degrees at the ends force an interesting planar substructure: An end of infinite vertex-degree produces the $\N\times\N$-grid as a minor (this is an old result of Halin~\cite{halin65}), and an end of vertex-degree at least $\frac 32k-1$ forces a $[k]\times \N$-grid-minor (and this bound is best possible). The latter result was not known before.

However, our notion of vertex-/edge-degrees is not strong enough to make  extensions of Theorems~\ref{thm:minor} and~\ref{thm:topminor}  possible. This can be seen by taking the infinite $r$-regular tree and inserting the edge set of some spanning subgraph at each level (Example~\ref{ex:Gk}). With a little more effort we can transform our
 example into one with infinitely many ends of large but finite vertex-/edge-degree (Example~\ref{ex:Gk'}).

To overcome this problem, we introduce in Section~\ref{sec:minors} a new end degree notion, the {\em relative degree}, that allows us to extend Theorems~\ref{thm:minor} and~\ref{thm:topminor} to infinite locally finite graphs  (Theorem~\ref{thm:minorLocFin}). Moreover, every locally finite graph of minimum degree/relative degree at least $k$ has a finite subgraph of average degree at least~$k$ (Theorem~\ref{thm:redu}).
An application of Theorem~\ref{thm:minorLocFin} is investigated in Section~\ref{sec:girth}, where we ask whether as in finite graphs, large girth can be used for forcing large complete minors.

\subsection{Minimal $k$-connectivity}

The subjects of the second part of our survey are minimally $k$-connected graphs. Minimality may here mean minimality with respect to either edge or vertex deletion, and it may also mean mimimality with respect to taking subgraphs. Minimality has been studied mainly for finite  graphs~\cite{lick,hamidoune,maderAtome,maderMin,maderEckenVom}, but also for infinite graphs~\cite{halinUnMin,maderUeberMin, minim}. See~\cite{minim} for an overview of results on edge-and vertex-minimality in finite graphs, see also~\cite{BBExtGT,frankHofC}.

Edge-minimally $k$-connected graphs, i.e.~those that are $k$-connected but lose this property upon the deletion of any edge, have received most attention in the literature and will be the subject of Section~\ref{sec:edgeDel}. It is known that these graphs have vertices of degree $k$. Even bounds on the number of such vertices are known~\cite{maderEckenVom,maderUeberMin}: Every {\em finite} edge-minimally $k$-connected graph must have at least $k+1$ vertices of degree $k$, and every  {\em infinite} edge-minimally $k$-connected graph $G$ has $|G|$ such vertices. Moreover, they appear on every (finite) cycle of~$G$. 

Unlike in finite graphs, however, infinite $k$-connected graphs need not have edge-minimally $k$-connected subgraphs. One example is the double-ladder (see the end of Section~\ref{sec:edgeDel}). 

It is thus natural to shift our investigations to certain `edge-minimally $k$-connected standard subspaces' which have the advantage that they do exist, at least in every locally finite $k$-connected graph (Lemma~\ref{lem:weak}). Then, most of the results for graphs mentioned above carry over to standard subspaces.
For the definition of these subspaces, we will have to we view the point set of a graph $G$ together with its ends as a {\em topological space}, see Section~\ref{sec:def}. 

Vertex-minimally $k$-connected graphs, i.e.~those graphs that are $k$-connected but lose this property upon the deletion of any vertex, are the topic of Section~\ref{sec:VDel}. It is known that finite such graphs have at least two vertices of `small' degree~\cite{hamidoune}, where `small'  now means $\frac 32k-1$ (which is best possible). This result carries over to infinite graphs, if we allow for  ends of small degree as well as vertices~\cite{minim}. It is necessary to allow also ends here.

Minimal $k$-connectivity with respect to taking subgraphs/induced subgraphs will be discussed in Section~\ref{sec:otherMin}. The respective questions for subspaces will be treated in Section~\ref{otherSubsp}.

Section~\ref{sec:EDel} and~\ref{sec:Vedge} investigate the same problems as Sections~\ref{sec:edgeDel} and~\ref{sec:VDel}, but for {\em edge-}connectivity. In Section~\ref{sec:EDel},  the existence and quantity of vertices of degree~$k$ in  edge-minimally $k$-edge-connected graphs (and sometimes multigraphs) are studied. In Section~\ref{sec:Vedge} we focus on vertex-minimally $k$-edge-connected graphs and multigraphs. The results shown in these two sections are taken from~\cite{minim}.

\subsection{Spanning circles and  trees}

The third and last part of the present survey, Section~\ref{sec:3}, deals with extremal problems concerning {\em circles}, {\em topological trees/forests}, and {\em arcs}, which shall be introduced in Section~\ref{sec:def}.  In addition to being natural extensions of the concepts of cycles, trees, forests, and paths in finite graphs, all these notions have proved over the last decade to be of immense use in infinite graph theory (see~\cite{DBook} or the survey~\cite{diestelBanffsurvey}).

Section~\ref{sec:Ham} presents results and open problems concerning {\em Hamilton circles}. The main result seems to be Georgakopolous' extension (Theorem~\ref{agelos}) of Fleischner's theorem that the square of any locally finite $2$-connected graph has a Hamilton cycle. The main conjecture, on the other hand, is due to Bruhn (Conjecture~\ref{henning}), and would extend a result of Tutte which states that every planar $4$-connected graph has a Hamilton cycle. For these and more problems/results, see Section~\ref{sec:Ham}.

In Section~\ref{sec:tree-p} we turn our attention to forests and spanning trees. These play the lead role in two well-known results from finite graph theory: the tree-packing theorem and the arboricity theorem. The former is about the number of edge-disjoint spanning trees of a graph. It states that if every partition of the vertex set of a graph $G$ is crossed by at least as many edges as $k$ edge-disjoint spanning trees would send across, then in fact, $G$ has $k$ edge-disjoint spanning trees. The infinite locally finite analogue is false for `traditional' spanning trees, but Bruhn and Diestel showed it holds true for topological spanning trees (Theorem~\ref{thm:locFinTreePack}). 

A related result, the arboricity theorem, 
 extends easily, if we do not require these forests to be topological ones, but extends also if we do (although then, a further condition is needed). See Section~\ref{sec:tree-p} for all details.

We close the last part of our survey with a topic that would have also fitted into the second part: connectivity-preserving arcs and cycles/circles. These are paths or cycles whose deletion does not reduce the connectivity `too much'.

A well known conjecture of Lov\'asz in this respect states that there is a function $f$ so that every finite $f(k)$-connected graph has an induced cycle so that the deletion of its vertices leaves the graph $k$-connected, and moreover, that one may prescribe an edge which the cycle has to contain. There are several weakenings and modifications of this conjecture which have been proved in finite graphs. We ask for extensions to infinite graphs in Section~\ref{sec:con-pres}.

\section{Terminology}\label{sec:def}

All our notation is as in~\cite{DBook}, but we take the oppotunity here to remind the reader of the few less standard concepts.

One of the main concepts in infinite graph theory is that of the {\em ends} of a graph $G$. An end of $G$ is an equivalence class of rays (i.e.~one-way infinite paths) of $G$, where we say that two rays are equivalent if no finite set of  vertices separates them. We denote the set of ends of a graph $G$ by $\Omega(G)$. 

The {\em vertex-degree} and the {\em edge-degree} of an end $\omega\in\Omega(G)$ were introduced in~\cite{degree} resp.~in~\cite{arbo}. Sometimes, one refers to both at the same time speaking informally of the {\em end degree}. 
The { vertex-degree} $d_v(\omega)$ of $\omega$ is defined as the maximum cardinality of a set of (vertex-)disjoint rays in $\omega$, and the { edge-degree} $d_e(\omega)$ of $\omega$ is defined as the maximum cardinality of a set of edge-disjoint rays in $\omega$. These maxima exist~\cite{halin65}, see also~\cite{DBook}. Cleary, the vertex-degree of an end is at most its edge-degree. We shall encounter a third end degree notion in Section~\ref{sec:minors}.

For a subgraph $H$ of a graph $G$, we write $\partial_v H := N(G-H)$ for its {\em vertex-boundary}. Similarly, $\partial_e H:= E(H,G-H)$ is the {\em edge-boundary} of $H$.

An induced connected subgraph $H$ of an infinite graph that has a finite vertex-boundary is called a {\em region}. If
$H$ contains rays of an end $\omega$, we will say that $H$ is a {\em region of} $\omega$.

For $k\in\N$, a separator of a graph of size $k$ will often be called a {\em $k$-separator}, and {\em $k$-cuts} are defined analogously. We say that a separator (or cut) $S$ of a graph $G$ separates some set $A\subseteq V(G)$ {\em from} an end $\omega\in\Omega(G)$, if the component of $G-S$ that contains rays of $\omega$ does not meet $A$.

\smallskip

The rest of this section is dedicated to the topological viewpoint on (infinite) graphs that has been introduced in~\cite{cyclesI,cyclesII,tst}. With a few exceptions we shall not use these concepts until Section~\ref{sec:subspace}, so the reader might wish to read the rest of this section only then.

We first define a topological space $|G|$ on the point set of the graph $G$ plus its ends. The topology is as on a $1$-complex, only that we allow as basic open neighbourhoods of a vertex $v$ only sets of half-open edges of the same length $\varepsilon$, and the basic open neighbourhoods of an end $\omega$ are defined  as follows. For each finite set $S\subseteq V(G)$, and for $\varepsilon>0$ let $C_{S,\omega}$ be the (unique) component of $G-S$ that contains rays of $\omega$. Let $\Omega_{S,\omega}$ be the set of all ends that have rays in $C_{S,\omega}$, and let $E_{S,\varepsilon,\omega}$ be the set of half-open intervals of length $\varepsilon$ of the edges in $E(S,V(C))$, one for each edge. Now $C_{S,\omega}\cup\Omega_{S,\omega}\cup E_{S,\varepsilon,\omega}$ is a basic open neighbourhood of $\omega$.

A {\em standard subspace} of $|G|$ is a closed subspace that contains every edge of which it contains an inner point. Observe that a standard subspace is thus nothing else than the closure $\overline H$ of a subgraph $H$ of $G$. Later in the text, we shall give a definition of end degrees in subspaces.

Finally we shall need the notion of circles and arcs, which are the infinite analogues of paths and cycles, and will be used mainly in Section~\ref{sec:3}. A {\em circle} is the homeomorphic image of the unit cycle in $|G|$ . An {\em arc} is the homeomorphic image of the unit interval. Observe that these definitions include the traditional cycles and paths. 

The definition of a circle gives rise to a new concept of trees and forests: These are now required to be void of circles (and not only finite cycles). Thus we define a {\em topological tree} in $G$ as a path-connected standard subspace of $|G|$ that contains no circles, and a {\em topological forest} as a union of such topological trees. 

\bigskip

\section{Degrees and substructure}

\subsection{Large complete minors in rayless graphs}\label{sec:rayless}

We start this section on substructures with an extension of Theorems~\ref{thm:minor} and~\ref{thm:topminor} to infinite rayless graphs. The functions $f_1$ and $f_2$ are as defined in these theorems.

\begin{theorem}\label{prop:rayless}
Let $G$ be a rayless graph. If each vertex of $G$ has degree at least $f_1(r)$, then $K^r$ is a minor of $G$, and if each vertex of $G$ has degree at least $f_2(r)$, then $K^r$ is a topological minor of $G$.
\end{theorem}

In fact, Theorem~\ref{prop:rayless} follows at once from Theorems~\ref{thm:minor} and~\ref{thm:topminor} together with the following reduction theorem:

\begin{theorem}\label{thm:reduRayless}
Let $G$ be a rayless graph of minimum degree $m$. Then $G$ has a non-empty finite subgraph of minimum degree $m$.
\end{theorem}

In order to prove Theorem~\ref{thm:reduRayless}, we need 
K\H onig's infinity lemma:
\begin{lemma}\label{inflemma}$\!\!${\bf\cite{DBook}}
 Let $G$ be a graph on the union of disjoint finite non-empty sets $S_i$, $i\in\N$, so that each $v\in S_i$ has a neighbour in $S_{i-1}$. Then $G$ has a ray.
\end{lemma}

 \begin{proof}[Proof of Theorem~\ref{thm:reduRayless}]
  We start with any finite non-empty vertex set $S_0$. For $i\geq 1$ we  choose for each vertex $v\in S_{i-1}$ a set $S_v$ of $\max\{0,m-d_{G[ \bigcup_{j< i} S_j]}(v)\}$ neighbours of $v$ in $V(G)\setminus \bigcup_{j< i} S_j$. This is possible, as by assumption $v$ has degree at least $m$ in $G$. We set $S_i:=\bigcup_{v\in S_{i-1}}S_v$. 
  
  Now if $S_i=\emptyset$ for some $i$, then $G[ \bigcup_{j< i} S_j]$ is the desired subgraph of $G$. On the other hand, if $S_i\neq\emptyset$ for all $i\in\N$,  we may apply Lemma~\ref{inflemma} to find a ray in $G$, a contradiction, as $G$ is rayless.
 \end{proof}

\subsection{Grid minors}\label{sec:grid}

From now on, we will deal with graphs that may have rays. We have already seen in the introduction that then large degrees at the vertices are not enough to force even cycles. We shall thus use additionally the end degrees in order to force interesting substructures in infinite graphs. In this subsection, we start modestly  by asking for minors that are planar. 

Particularly interesting planar graphs are the grids. The {\em infinite grid} $\mathbb Z\times \mathbb Z$ is the graph on $\mathbb Z^2$ having all edges of the form $(m,n)(m+1,n)$ and of the form $(m,n)(m,n+1)$, for $m,n\in\mathbb Z$. The {\em half-grid} $\N\times \Z$, the  {\em quarter-grid} $\N\times \N$, and the $[k]\times \N$-grid are the induced subgraphs of $\Z\times\Z$ on the respective sets.

 A well-known result in infinite graph theory concerns the quarter-grid\footnote{Observe that when considering minors, it makes no difference whether we work with the half-grid or the quarter-grid, since, as one easily checks, each of the two is a minor of the other.}, which is a minor of every graph that has an end of infinite vertex-degree (this is a classical result of Halin~\cite{halin65} who called such ends {\em thick ends}).

\begin{theorem}[Halin~\cite{halin65}]
Let $G$ be graph which has an end $\omega$ of infinite vertex-degree. Then the $\mathbb N\times \mathbb N$-grid is a minor of $G$.
\end{theorem}

From Halin's proof it follows that the rays of the subgraph of $G$ that can be contracted to $\mathbb N\times\N$  belong to $\omega$ (see also the proof in Diestel's book~\cite{DBook}). On the other hand, it is clear that if a subdivision of the  quarter-grid appears as a subgraph of some graph $G$, the its rays belong to an end of infinite vertex-degree in $G$.

Thus, it is not surprising that assuming large (but not infinite) degrees and vertex-degrees we cannot force a quarter-grid minor.
One example for this fact is $\tilde G_k$ which is to be defined after Theorem~\ref{thm:hcs}, another, even planar, example is the graph $G_k'$ from Example~\ref{ex:Gk'}.

However, both graphs  contain something quite similar to a quarter-grid: a $[k]\times \mathbb N$ grid, where $k$ depends on the minimum vertex-degree we required at the ends. In fact, such a grid always appears in a graph with an end $\omega$ of large enough vertex-degree. It will follow from the proof that the rays corresponding to the rays of the minor, in $G$ belong to $\omega$.

\begin{theorem}\label{lem:rays}
Let $k\in\N$ and let $G$ be graph which has an end $\omega$ of vertex-degree at least~$\frac 32 k-1$. Then the $[k]\times \mathbb N$-grid is a minor of $G$.
\end{theorem}

The bound on the vertex-degree is sharp. This is illustrated by Example~\ref{ex:nogrid}, after the proof of Theorem~\ref{lem:rays}.

\begin{proof}[Proof of Theorem~\ref{lem:rays}]
We shall proceed by induction on $k$. For $k=1$ and $k=2$, the assertion clearly holds, so assume that $k\geq 3$ and that $\omega$ is an end of a graph $G$ with $d_v(\omega)\geq \frac 32 k-1$.

 Choose a set $\mathcal R$ of $d_v(\omega)$ disjoint rays from $\omega$.  Consider the auxiliary graph $H$ with $V(H):=\mathcal R$ where two vertices $R$ and $R'$ are adjacent if there exists an infinite set of disjoint $V(R)$--$V(R')$ paths in $G$ which avoid all $R''\in\mathcal R$ with $R''\neq R,R'$. Let $T$ be a spanning tree of $H$. Clearly, if $T$ happens to be a path, it is easy to construct the desired minor. 

So suppose otherwise. Then $T$ has (at least) three leaves $R_1$, $R_2$, $R_3$. Observe that the graph $G':=G-V(\bigcup_{j=1,2,3} R_j)$  has an end $\omega '$ of degree 
\[
 d_v(\omega ')=d_v(\omega)-3\geq\frac 32 k -4=\frac 32(k-2)-1
\]
whose rays, when viewed in $G$, belong to $\omega$.
Hence, by induction, the $[k-2]\times\N$-grid is a minor of $G'$. In other words, $G'$ contains a set of rays $Q_1,Q_2,\ldots Q_{k-2}\in\omega'$, and furthermore, each $Q_i$ is linked to $Q_{i+1}$ by infinitely many disjoint paths, which do not meet any other $Q_j$.

In $G$, the $Q_i$ belong to $\omega$. Thus, since $|\mathcal R|=d_v(\omega)$, each $Q_i$ meets $\bigcup\mathcal R$ infinitely often. Hence each $Q_i$ meets (at least) one of the rays in $\mathcal R$, which we shall denote by $R(Q_i)$, infinitely often.

The tree $T$ from above contains three paths $P_i$, $i=1,2,3$, so that $P_i$ starts in $V(R_i)$ and ends in $\bigcup_{i=1}^{k-2}V(R(Q_i))$. Since $R_1$, $R_2$ and $R_3$ are leaves of $T$, the $P_i$ can be chosen so that they are disjoint except possibly in their endvertices. Using the path systems in $G$ represented by the $P_i$, it is now easy to see that for each $R_j$, $j=1,2,3$, there is a $Q_{i_j}$ among the $Q_i$ such that there exist an infinite family of disjoint $V(R_j)$--$V(Q_{i_j})$ paths which avoid all other $Q_{i'}$ and $R_{j'}$. Say $i_1\leq i_2\leq i_3$.

In order to see that the  $[k]\times \mathbb N$-grid is a minor of $G$, we shall now define a family of rays $\tilde Q_1,\tilde Q_2,\ldots \tilde Q_{k}\in\omega$ so that $\tilde Q_i$ and $\tilde Q_{i+1}$ are connected by infinitely many disjoint paths which do not meet any other $\tilde Q_i$. For $i< i_1$ set $\tilde Q_i:=Q_i$, and for $i>i_3+2$ set $\tilde Q_{i}:=Q_{i-2}$. Set $\tilde Q_{i_2+1}:=R_2$. For $i\neq i_2+1$ with $i_1<i< i_3+2$, we choose $\tilde Q_i$ as a suitable ray which alternatively visits $Q_{i-1}$ and $Q_i$, if $i\leq i_2$, or $Q_{i-2}$ and $Q_{i-1}$, if $i>i_2+1$. Finally, $\tilde Q_{i_1}$ and $\tilde Q_{i_3+2}$ are chosen so that they alternate between $R_1$ and $Q_{i_1}$, respectively between $Q_{i_3}$ and $R_3$. Clearly this choice of the rays $\tilde Q_i$ ensures that, together with suitable connecting paths, the $\tilde Q_i$ may be contracted to a $[k]\times \mathbb N$-grid.
\end{proof}

\begin{figure}[ht]
      \centering
      \includegraphics[scale=0.65]{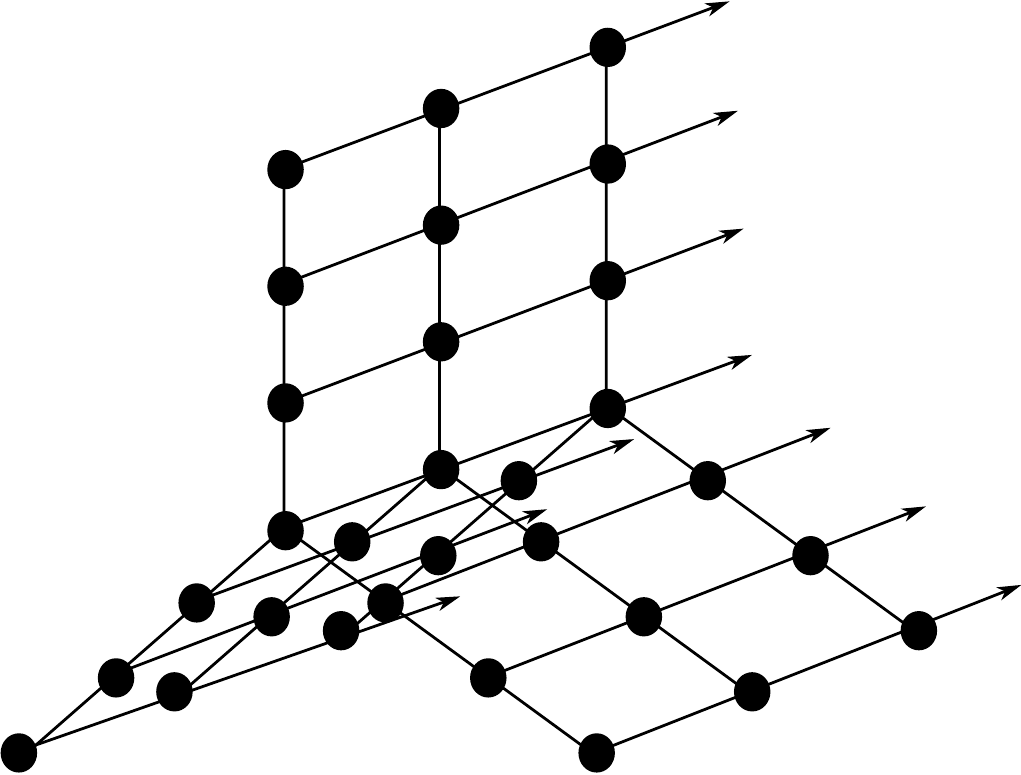}
      \caption{The graph $Y(3)$ from Example~\ref{ex:nogrid}.\label{fig:Y}}
\end{figure}

\begin{example}\label{ex:nogrid}
Denote by $K_{1,3}(\ell)$ the graph that is obtained by replacing each edge of $K_{1,3}$ with a path of length $\ell$. Define $Y(\ell):=K_{1,3}(\ell)\times \N$. (That is, for each $i\in\N$, we take a copy of $K_{1,3}(\ell)$ and add an edge between every $i$th and $(i+1)$th copy of each vertex in $K_{1,3}(\ell)$.)\\
Clearly, the vertex-degree of the unique end of $Y$ is $3\ell +1$. We shall show in Lemma~\ref{lem:nogrid} that the $[k]\times \mathbb N$-grid is not a  minor of $Y(\ell)$, for $k=2\ell +2$. 
\end{example}

\begin{lemma}\label{lem:nogrid}
Let $\ell\in\N$ and let  $k=2\ell +2$.
Then the graph $Y(\ell)$ from Example~\ref{ex:nogrid} has an end of vertex-degree $\frac32k-2$, but  the $[k]\times \mathbb N$-grid is not a  minor of $Y(\ell)$.
\end{lemma}

\begin{proof}
Suppose otherwise. Then the graph $Y(\ell)$ contains a family of rays $\mathcal R:=\{R_1, R_2,\ldots R_k\}$ such that for $i=1,2,\ldots, k-1$, there are infinitely many finite paths connecting  $R_i$ with $R_{i+1}$, such that all these paths are all disjoint, except possibly in their endvertices, and such that they avoid all $R_{i'}$ with $i'\neq i,i+1$.\\
 Let $n\in\N$ be such that all $R_i$ meet $Y_n:=K_{1,3}(\ell) \times \{ n\}$, the $n$th copy of $K_{1,3}(\ell)$ in $Y(\ell)$. Write  $V(Y_n)$ as $\{v_0, v^1_1, v^1_2,\ldots, v^1_\ell, v^2_1, v^2_2,\ldots, v^2_\ell, v^3_1, v^3_2,\ldots, v^3_\ell\}$ where each $v_0v^j_1v^j_2\ldots v^j_\ell$ induces a path in $Y_n$.  \\
For each $j=1,2,3$ consider that ray $R(j)\in\mathcal R$ that meets a $v^j_m$ with largest index $m$. Observe that (at least) one of these three rays, say $R(1)$ is neither equal to $R_1$ nor to $R_k$. Let $R'(1)$ be the ray in $\mathcal R$ that  meets $v^1_m$ with the second largest index $m$, or, if there is no such, let $R'(1)$ be the ray that meets~$v_0$ (which then exists, since $|\mathcal R|=k>2\ell+1$ and since each ray of $\mathcal R$ meets $Y_n$).\\
 We claim that $S:=V(R'(1))\cup V(\bigcup_{h\leq n} Y_h)$ separates $R(1)$ from the rest of the $R_i$, which clearly leads to the desired contradiction, since $R(1)\neq R_1,R_k$, and thus has to be connected to two of the $R_i$ by infinitely many disjoint finite paths that avoid all other $R_i$. So suppose otherwise, and let $P$ be a path that connects $R(1)$ in $Y(\ell)-S$ with some $R_{i^*}\in\mathcal R$.\\
By construction of $Y(\ell)$, this is only possible if $R'(1)$ uses vertices of the type $v^2_m$ or $v^3_m$. Let $\tilde n$ be the smallest index $\geq n$ such that this occurs, say the $\tilde n$th copy of $v^2_ 1$ lies on $R'(1)$. Then also the $\tilde n$th copy of $v_0$ lies on $R'(1)$, and furthermore,  all other $R_i$ (with the exception of $R(1)$) have to pass through the $\tilde n$th copies of the vertices $v^2_2, v^2_3,\ldots, v^2_\ell, v^3_1, v^3_2,\ldots, v^3_\ell$. Hence the total number of rays in $\mathcal R$ cannot exceed $2\ell +1$, a contradiction, as $k=2\ell +2$.
\end{proof}

\subsection{Highly connected subgraphs}\label{sec:hcs}

We shall now see another example of how large end degrees and large degree at the vertices force a certain dense substructure. In fact, assuming large degree and large vertex-/edge-degree we can ensure highly connected or highly edge-connected subgraphs in infinite graphs. This is the content of Theorem~\ref{thm:hcs} below. Before we state it, let us quickly remark that conversely, in a locally finite 
$k$-connected/$k$-edge-connected graph, all ends have vertex-/edge-degree at least $k$. 
This will follow at once from Lemma~\ref{lem:converse} of Section~\ref{sec:Vtrans}.

\begin{theorem}\label{thm:hcs}$\!\!${\bf\cite[{\rm Theorems 3.1 and 5.1}]{hcs} }
Let $k\in\mathbb N$, let $G$ be an infinite graph.
\begin{enumerate}[(a)]
\item If all vertices of $G$ have degree greater than $2k(k+1)$, and all ends of $G$ have vertex-degree at least $2k(k+3)$, then $G$ contains  a $(k+1)$-connected subgraph.
\item If all vertices of $G$ have degree at least $2k$, and all ends of $G$ have edge-degree at least $2k$, then $G$ contains  a $(k+1)$-edge-connected subgraph.
\end{enumerate}
\end{theorem}

We can find these highly (edge-)connected subgraphs inside every region of $G$~\cite{hcs}.

Observe that we have no control on whether the highly connected subgraph from Theorem~\ref{thm:hcs} is finite or infinite. That is, there are examples of graphs with the prescribed degrees but no finite highly connected subgraph, and others which have no infinite highly connected subgraph. An example of the first kind is the graph $G_k$ from Example~\ref{ex:Gk}.
An example of the second kind can be constructed as follows. 
For  $k\in\mathbb N$, consider the $k\times \mathbb N$ grid, and for each vertex $v$ of this grid, take a copy of $K^{k+1}$, and identify one of its vertices with $v$. The obtained graph $\tilde G_k$ has an end of vertex- and edge-degree $k$, and all vertices have degree at least $k$, but $\tilde G_k$ has no infinite $5$-connected subgraph. Similar examples can be constructed for edge-connectivity.

There are also examples of graphs (for all $k\in\mathbb N$) that have very large degree at the vertices, and a degree of order $k\log k$ at the ends, but no $k$-connected subgraphs~\cite{hcs}. Thus the at first sight surprising quadratic bound on the vertex degrees of the ends is not too far from best possible.

\begin{theorem}$\!\!${\bf\cite[{\rm Theorem 6.1}]{hcs} }
 For each $k=5\ell$, where $\ell\in\N$ is even, there exists a locally finite
graph whose vertices have degree at least $2^\ell$, whose ends have vertex-degree at
least $\ell \log \ell$, and which has no $(k + 1)$-connected subgraph.
\end{theorem}

So, Theorem~\ref{thm:hcs} can not be improved in this sense. We may ask however, whether the theorem holds for standard subspaces of infinite graphs. The analogue of this question in finite graphs would be to ask whether the theorem stays true for subgraphs, which is obviously true. We see that the infinite setting allows for more subtleties than the finite one.

For brevity, we shall only concentrate on part a) of Theorem~\ref{thm:hcs}, that is, the vertex-version. First, we have to define a notion of $k$-connectivity for standard subspaces. 
There are two options which seem natural. 

Call a path-connected standard subspace $X$ of the space $|G|$ (that is associated to some graph $G$) that contains at least $k+1$ vertices of $G$ {\em strongly $k$-connected} if deleting up to $k$ ends or vertices (the latter together with all adjacent edges) leaves $X$ path-connected. Call $X$ {\em weakly $k$-connected}, if deleting up to $k$ vertices together with all adjacent edges from $X$ leaves a path-connected space. Observe that for $X=|G|$, our notions coincide, and coincide with $k$-connectivity of $G$.

We also have to define the degree of an end $\omega$ in the standard subspace $X$. Following~\cite{degree}, we say that the {\em vertex-degree} of $\omega$ in $X$ is the maximum of the cardinalities of the sets of arcs
in $X$ that are disjoint except in their common endpoint $\omega$. The {\em edge-degree} in $X$ is defined analogously.

Then, the question is whether there exists a function $f:\mathbb N\to \mathbb N$ so that if $Y$ is a standard subspace of an infinite graph $G$ whose vertices and ends all have degree resp.~vertex-degree at least $f(k)$ in $Y$, then there is a standard subspace $X\subseteq Y$ which is weakly or even strongly $k$-connected.

For strong $k$-connectivity, the answer is no. This is illustrated by the following simple example.

\begin{example}\label{ex:hcsSubsp}
Consider the graph $G_k$ which we obtain from the $k$-regular tree $T_k$ by inserting a spanning cycle at each level. We consider the standard subspace $\overline T_k$ of $|G_k|$. \\ Clearly, all vertices and also the unique end $\omega$ of $\overline T_k$ have (vertex-)degree at least $k$. However, for any standard subspace $X$ of $\overline T_k$ that contains at least $3$ vertices, we can choose one of its vertices so that its deletion  plus the deletion of $\omega$ (if present in $X$) destroys the path-connectivity.
\end{example}

However,  Theorem~\ref{thm:hcs} might still extend to standard subspaces if we use the weaker notion of $k$-connectivity:

\begin{problem}
 Is there a function $f:\mathbb N\to \mathbb N$ so that the following holds: If $Y$ is a standard subspace of an infinite graph $G$, such that all vertices and ends of $Y$ have degree resp.~vertex-degree at least $f(k)$ in $Y$, then $G$ has a  weakly $k$-connected standard subspace $X\subseteq Y$?
\end{problem}

\subsection{Connectivity of vertex-transitive graphs}\label{sec:Vtrans}

Let us now pose the question from the previous section for vertex-transitive graphs. As vertex-transitive graphs are regular, we need no longer use the term `mimimum degree'. Thus our question from Section~\ref{sec:hcs}  reduces to the following in vertex-transitive graphs: Which degree at each vertex do we need in order to ensure that our graph has a $k$-(edge-)connected subgraph?

It is known that in finite graphs a degree of $k$ is enough, and moreover the subgraph will be the graph itself. In fact, every finite vertex-transitive $k$-regular connected graph is $k$-edge-connected~\cite{maderMin}. It is even $k$-connected, as long as it does not contain $K^4$ as a subgraph~\cite{maderSym}.

In infinite graphs, this is no longer true, if we only require degree $k$ at the vertices, because of the trees. However, if we require a vertex-/edge-degree of at least $k$ at the ends (which is conversely implied by the $k$-(edge-)connectivity, see below), we can obtain analogous results for infinite locally finite graphs. We may even drop the condition on the degrees of the vertices.

\begin{proposition}\label{prop:Vtrans}
Let $G$ be an infinite locally finite graph, let  $k\in\mathbb N$. Suppose that $G$ is vertex-transitive and connected.
\begin{enumerate}[(a)]
\item $G$ is $k$-connected if and only if all ends of $G$ have vertex-degree at least $k$.
\item $G$ is $k$-edge-connected if and only if all ends of $G$ have edge-degree at least~$k$.
\end{enumerate}
\end{proposition}

In fact, the forward  implications in Proposition~\ref{prop:Vtrans} are easily implied by the following result, whose proof is not very difficult and can be found in~\cite{degree} for the edge-case (the vertex-case is analogous). 

\begin{lemma}\label{lem:converse}
 Let $k\in \N$, let $G$ be a locally finite graph, and let $\omega\in\Omega(G)$. Then 
\begin{enumerate}[(i)]
\item $d_v(\omega)=k$ if and only if $k$ is the smallest integer such that every finite set $S \subseteq V (G)$
can be separated from $\omega$ with a $k$-separator, and
\item $d_e(\omega)=k$ if and only if $k$ is the smallest integer such that every finite set $S \subseteq V (G)$
can be separated from $\omega$ with a $k$-cut.
\end{enumerate}
\end{lemma}

\begin{proof}[Proof of Proposition~\ref{prop:Vtrans}]
Because of Lemma~\ref{lem:converse} we only need to prove the backward implications. Let us only prove the implication for (a), for (b) this is analoguous. 

Suppose the implication is not true, and let $S$ be an $\ell$-separator of $G$, for some $\ell< k$. Choose a vertex $w$ at distance at least $\max\{dist(u,v):u,v\in S\}+1$ from all $v\in S$. (Observe that such a vertex $w$ exists, since $G$ is infinite, locally finite and connected.) 
Now, let $\phi$ be an automorphism of $G$ that maps some vertex from $S$ to $w$. Then $\phi(S)$ is contained in one component of $G-S$. 

Next, choose an automorphism $\phi'$ that maps $\phi (w)$ `far away' from $\phi (S)$ to a component of $G-\phi(S)$ that does not contain $S$. Continuing in this manner, we arrive at a sequence  $S$, $\phi(S)$, $\phi'(\phi((S))$, $\phi''(\phi'(\phi(S)))$, $\ldots$ of $\ell$-separators of $G$.
It is not difficult to construct a ray that meets each of these separators and hence defines an end of vertex-degree $\ell <k$. This contradicts our assumption that all ends have vertex-degree at least $k$.
\end{proof}

\subsection{Two counterexamples}

This short section is dedicated to two examples which show that large degree and large vertex-degree together are not strong enough assumptions to force large complete minors. The difference between the two examples is that the latter does not have ends of infinite vertex-degree.

\begin{example}\label{ex:Gk}
For given $k$, take the $k$-regular tree $T_k$ with levels $L_0,  L_1,L_2,\ldots$ and insert the edge set of a spanning cycle $C_i$ at each level $L_i$ of $T_k$  (cf.~Example~\ref{ex:hcsSubsp}). This can be done in a way so that the obtained graph $G_k$ is still planar. \\
Clearly, $G_k$ has one end of infinite vertex- and edge-degree, and furthermore, all vertices of $G_k$ have degree at least $k$.
It is easy to see that $G_k$ is $k$-connected, but being planar, $G_k$ has no complete minor of order greater than $4$.
\end{example}
%

By deleting some (carefully chosen) edges from $G_k$, we obtain a planar graph of high minimal degree and vertex-degree whose (continuum many) ends all have finite vertex-degree:

\begin{example}\label{ex:Gk'}
Let $k\in \mathbb N$ be given, and consider 
the graph $G_k$ from Example~\ref{ex:Gk}. Now, for each $i\in\mathbb N$, delete the edge $vw\in E(C_i)$ from $E(G_k)$, if $v$ and $w$ have no common ancestors in the levels $L_{i-k+2}, L_{i-k+3},\ldots,L_{i-1}$. Denote the obtained graph by $G'_k$. \\
As $G_k$ is planar, also $G_k'$ is. Clearly, $k\leq d(v)\leq k+2$ for each $v\in V(G)$. We show in Lemma~\ref{lem:Gk'} that the ends of $G_k'$ have large, but finite vertex-degree.
\end{example}

\begin{figure}[ht]
      \centering
      \includegraphics[scale=0.35]{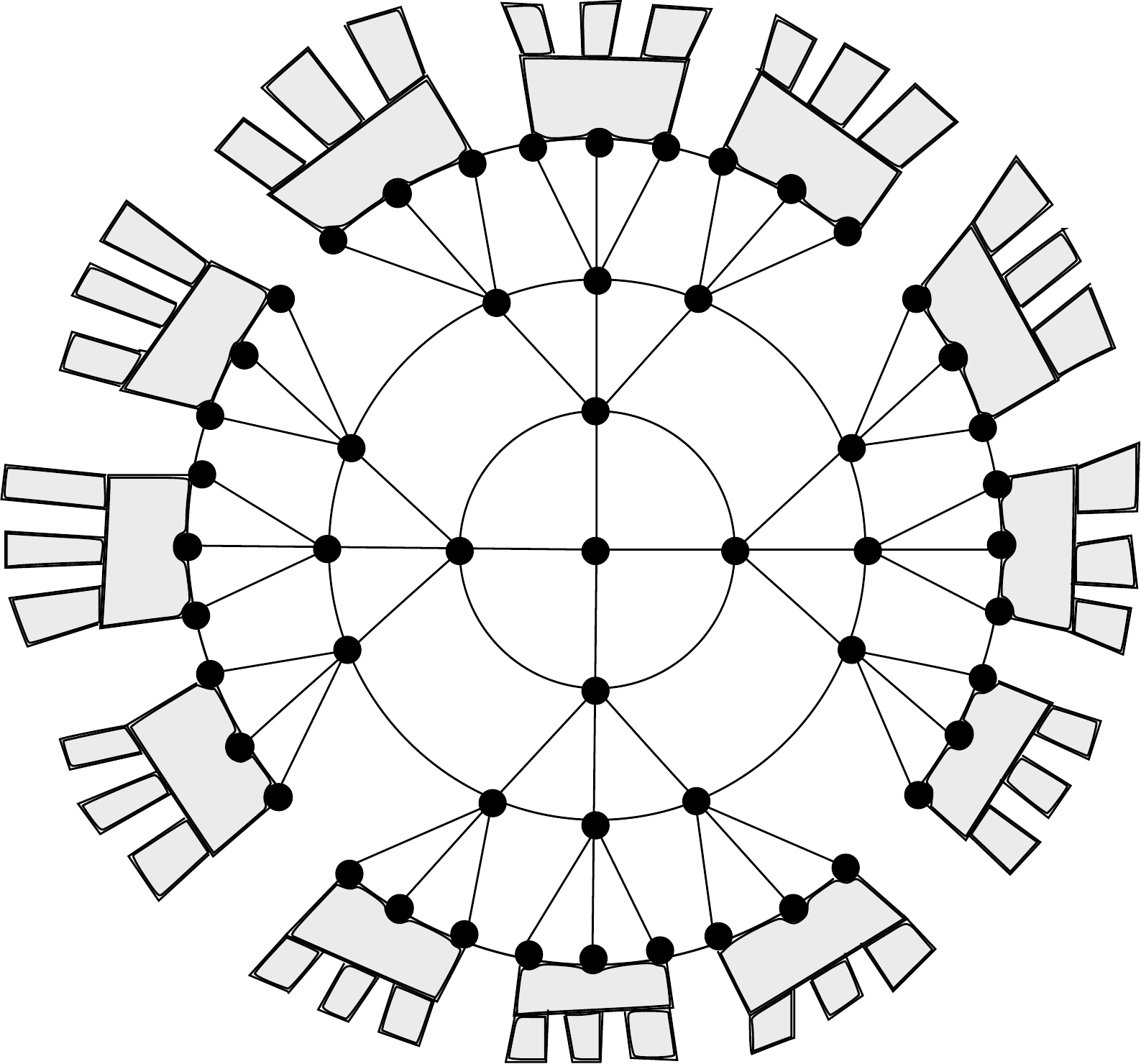}
      \caption{The graph $G_k'$ from Example~\ref{ex:Gk'} for $k=4$.\label{fig:Gk'}}
\end{figure}

\begin{lemma}\label{lem:Gk'}
The ends of the graph $G'_k$ from Example~\ref{ex:Gk'} all have vertex-degree between $k-2$ and $2k-3$.
\end{lemma}

\begin{proof}
Consider, for each $x\in V(G_k')$ the set
\[
 S_x:=\{x\}\cup\bigcup_{i=1,\ldots,k-2}N^i(x),
\]
where $N^i (x)$ here denotes the $i$th neighbourhood of $x$ in level $L_{m+i}$, supposing that $x$ lies in the $m$th level (of $T_k$).\\
Clearly for each $x\in V(G'_k)$, the set $S_x$ separates $G'_k$. Hence, already $\partial_v S_x$, which has order between $k-1$ and $2k-3$, separates $G'_k$.\\
Let us use the sets $S_x$ in order to show that the ends of $G'_k$ correspond to the ends of $T_k$. In fact, all we have to show is that for each ray $R\in G'_k$ there is a ray $R_T$ in $T_k$ that is equivalent to $R$ in $G_k'$. We can find such a ray $R_T$ by considering for each $i$ large enough the last vertex $v_i$ of $R$ in  $V(L_i)$. Now, $v_i\in S_{w_i}$ for exactly one $w_i\in V(L_{i-k+2})$. By definition of the $v_i$, the $w_i$ are adjacent to their successors $w_{i+1}\in V(L_{i-k+3})$. So,  $R_T:=w_kw_{k+1}w_{k+2}\ldots$ is a ray in $T_k$ as desired.

Thus $G'_k$ has continuum many ends, all of which have vertex-degree at most $2k-3$, because of the separators $\partial_v S_x$. 
It remains to show that each end $\omega$ of $G'_k$ has vertex-degree at least $k-2$. 

For this, fix $\omega\in\Omega(G)$ and  consider the union $S_\omega$ of the sets $\partial_v S_{w_i}$ for the ray $R=w_0w_1w_2w_3\ldots$ of $T_k$ that lies in $\omega$, where we assume that $R$ starts in $L_0=\{w_0\}$. By Lemma~\ref{lem:converse}, in order to see that $\omega$ has vertex-degree at least $k-2$ in $G[S]$ (and thus in $G$) we only have to show that no set of less than $k-2$ vertices separates $L_0$ from $\omega$ in $G[S]$.

So suppose otherwise, and let $T$ be such a separator. Since every vertex of $S$ has at least $k-1$ neighbours in the next level, we can reach the $2$nd, $3$rd, \dots ${k-2}$th level from $w_0$ in $G[S]-T$. By definition of $G'_k$, these levels contain spanning cycles, and thus, as $|T|>k-2$, there is a $w_i$ with $i\in\{1,2,\ldots,k-2\}$ which can be reached from $w_0$ in $G[S]-T$. We repeat the argument with $w_i$ in the role of $w_0$, observing that in $G[S]\cap (L_{i+1}\cup L_{i+2}\cup\ldots\cup L_{i+k-2})$, each level contains spanning paths, by construction of $G'_k$.
\end{proof}

\subsection{Large relative degree forces large complete minors}\label{sec:minors}

In the previous sections we explored which substructures may or may not be forced in an infinite graph if we assume large (vertex-)degree at both vertices and ends. In particular we saw that Theorems~\ref{thm:minor} and~\ref{thm:topminor} (with the average degree replaced by the minimum degree) do not extend  to infinite graphs that have rays. 

In the present section we shall overcome this problem. We will see that with a different, more appropriate notion of the end degree a satisfactory extension of Theorems~\ref{thm:minor} and~\ref{thm:topminor} to locally finite graphs is possible.

For this, let us first take a closer look at the graph $G_k'$ from Example~\ref{ex:Gk'}. Why do the large (vertex-)degrees not interfere with the planarity?
Observe that, for each finite set $S\subseteq V(G)$, the  edge-boundary of the subgraph $G_k'-S$ has about the same size as its vertex-boundary. So locally the density is never large enough to force non-planarity. Similar as in the tree $T_k$, the density that the high degrees should generate gets lost towards infinity. 

In order to avoid this behaviour, we have to prohibit regions $R$ of an end $\omega$ which have the property that $|\partial_e R|/|\partial_v R|$ is small, or at least we should prohibit sequences of such regions `converging' to $\omega$. This motivates us to define the {\em relative degree} of an end as the limit of the ratios above for `converging' sequences. This is not unnatural: applied to vertices this gives the usual degree, as each vertex $v$ is contained in a smallest region, namely $R=\{ v\}$, for which $|\partial_e R| /|\partial_v R|=d(v)$.

Let us make our idea more precise. Suppose that $G$ is a locally finite graph. We introduce a useful notation: if $H$ is a region of $G$, let us write $\Omega^G(H)$ for the set of all ends of $G$ that have rays in $H$.

 Now, write $(H_i)_{i\in\N}\rightsquigarrow\omega$ if $(H_i)_{i\in\N}$ is an infinite sequence of distinct regions of $G$ with $H_{i+1}\subseteq H_i-\partial_v H_i$ such that $\omega\in\overline H_i$  for each $i\in\N$. If moreover, $\partial_v H_{i+1}$ is an $\subseteq$-minimal $\partial_vH_i$--$\Omega^G (H_{i+1})$ separator for each $i\in\N$, then we write $(H_i)_{i\in\N}\rightarrow\omega$. Observe that such sequences always exist, as $G$ is locally finite. Define the {\em relative degree} of an end as
\[
\df(\omega):=\inf_{(H_i)_{i\in\N}\rightarrow\omega}\liminf_{i\rightarrow\infty}\frac{|\partial_e H_i|}{|\partial_v H_i|}.
\]
Note that it does not matter whether we consider the $\liminf$ or the $\limsup$, because if $(H_i)_{i\in\N}\rightarrow\omega$, also  all subsequences of $(H_i)_{i\in\N}$ converge to $\omega$. For the same reason we could restrict our attention to sequences $(H_i)$ where $\lim_{i\rightarrow\infty}\frac{|\partial_e H_i|}{|\partial_v H_i|}$ exists.

We remark  that if in the definition of the relative degree we replaced $(H_i)_{i\in\N}\rightarrow\omega$  with $(H_i)_{i\in\N}\rightsquigarrow\omega$, then the result would be a `degree' of $1$ for every end in any graph. Indeed, let  $(H_i)_{i\in\N}$ with $(H_i)_{i\in\N}\rightsquigarrow\omega$, and let $v_i\in\partial_v H_{3i}$ for $i\in\N$. Then the $v_i$ do not have common neighbours. We construct a sequence $(H'_j)_{j\in\N}$ with $H'_0:=H_0$, and, for $j>0$, we let $H'_j:=H_{i_j}-V_j$ where $i_j$ is such\footnote{For instance set $i_j:=\max\{dist(v,w)|v\in\partial_vH_0, w\in\partial_vH'_{j}\}+1$.} that $H_{i_j}\subseteq H'_j-\partial_v H'_j$, and $V_j$ consists of $j|\partial_e H_{i_j}|$ vertices $v_i$ with $i\geq i_j$. Then $(H'_j)_{j\in\N}\rightsquigarrow\omega$, and 
\[
 \liminf_{j\rightarrow\infty}\frac{|\partial_e H'_j|}{|\partial_v H'_j|}=\liminf_{j\rightarrow\infty}\frac{|\partial_e H_{i_j}|+\sum_{v\in V_j}d(v)}{|\partial_v H_{i_j}|+\sum_{v\in V_j}d(v)}=1.
\]
This shows that the additional condition that $\partial_v H_{i+1}$ is an $\subseteq$-minimal $\partial_vH_i$--$\Omega^G (H_{i+1})$ separator is indeed neccessary for the relative degree to make sense. For more discussion of our notion, see~\cite{LCM}.

Note that by Lemma~\ref{lem:converse}, in locally finite graphs, we can also express our earlier notions, the vertex- and the edge-degree, using converging sequences of regions. Here the $\rightsquigarrow$-convergence suffices:

\begin{align*}
 d_v(\omega)&=\inf_{(H_i)_{i\in\N}\rightsquigarrow\omega}\liminf_{i\rightarrow\infty}{|\partial_v H_i|}, \\
 d_e(\omega)&=\inf_{(H_i)_{i\in\N}\rightsquigarrow\omega}\liminf_{i\rightarrow\infty}{|\partial_e H_i|}.\label{d_e}
\end{align*}

Note that while $d_e\geq d_v$, there is no relation between $\df$ and any of $d_e$, $d_v$. Examples are not difficult to construct. For instance, take the union of complete graphs on $k$ vertices, one for each $i\in\N$, that gives a graph $H$ with vertex set $\bigcup_{i\in\N}\bigcup_{j=1}^k\{v_j^i\}$. Adding all edges $v_1^iv_1^{i+1}$ and $v_1^iv_2^{i+1}$, we obtain a graph with an end of vertex-degree $1$, edge-degree $2$ and relative degree $k$. On the other hand, adding to $H$ the edges $v_j^iv_j^{i+1}$ for all $i\in\N$ and all $j=1,\ldots ,k$, we get a graph with an end of vertex-/edge-degree $k$ and relative degree $1$. See~\cite{LCM}.

With the notion of the relative degree at hand, we can prove a very useful reduction theorem:

\begin{theorem}\label{thm:redu}
Let $G$ be a locally finite graph such that each vertex has degree at least $k$, and for each end $\omega$ we have $\df(\omega)\geq k$. Then~$G$ has a finite subgraph $H$ of average degree at least~$k$.
\end{theorem}

\begin{proof}
Choose a vertex $v\in V(G)$ and set $S_0:=\{v\}$. Inductively we shall construct a sequence $(S_i)_{i\in\mathbb N}$ of finite vertex sets with $S_i\subseteq S_{i+1}$ for all $i\in\N$. In each step $i\geq 0$ we start by considering the set $\mathcal A_i$ of all components $A$ of $G-(S_i\cup N_G(S_i))$. Let $\mathcal B_i\subseteq \mathcal A_i$ be the set of all those $B\in\mathcal A_i$ that contain a ray. Observe that as $G$ is locally finite,  each $A\in\mathcal A_i\setminus \mathcal B_i$ is finite. Moreover, since we may assume $G$ to be connected, $|\mathcal A_i\setminus \mathcal B_i|<\infty$, and thus $F_i:=\bigcup( \mathcal A_i\setminus\mathcal B_i)$ is finite.

Next, let $\mathcal C_i$ be the set of all components of $\bigcup\mathcal B_i\cup N_G(\bigcup\mathcal B_i)$. Note that $\mathcal C_i$ is finite and that for each $C\in \mathcal C_i$ 
\begin{equation}\label{minsepp}
\partial_v C\subseteq N_G(S_i)\text{ is an $\subseteq$-minimal $N_{S_i}(C)$--$\Omega^G (C)$ separator.}
\end{equation}
Let $\mathcal D_i\subseteq \mathcal C_i$ be the set of all those $D\in\mathcal C_i$ with
\begin{equation*}\label{ev<k}
\frac{|\partial_e D|}{|\partial _v D|}<k.
\end{equation*}
 Finally, set
\[
S_{i+1}:=S_i\cup F_i\cup (N_G(S_i)\setminus \bigcup\mathcal C_i)\cup\bigcup_{D\in\mathcal D_i}\partial_v D.
\]
This finishes the definition of the sets $S_i$. Note that by construction, $\partial_vC\subseteq N_G(S_i)\cap N_G(S_{i+1})$ and $N_G(C)\subseteq S_{i+1}\setminus \bigcup \mathcal D_i$ for each $C\in\mathcal C_i\setminus \mathcal D_i$. Hence, it is easy to show by induction that 
\begin{equation}\label{CsmD}
\mathcal C_i\setminus \mathcal D_i\subseteq \mathcal C_{i+1}\setminus \mathcal D_{i+1}\text{ for all }i\in\N.
\end{equation} 

Now, if there is an $i\in\N$ so that $\mathcal D_i=\emptyset$, then  $H:=G[S_{i+1}\cup N_G(S_{i+1})]$ is as desired. Indeed, then $N_G(S_{i+1})=\bigcup_{C\in\mathcal C_i\setminus \mathcal D_i}\partial_v C$. Thus by construction, and  by definition of $\mathcal D_i$, $H$ has average degree $\geq k$.

Otherwise, that is, if $\mathcal D_i\neq\emptyset$ for all $i$, we apply K\H onig's infinity lemma (Lemma~\ref{inflemma}) to the graph with vertex set $\bigcup_{i\in\N}\mathcal D_i$ which has an edge $CD$ whenever $C\in\mathcal D_i$, $D\in\mathcal D_{i+1}$ and $D\subseteq C$. Note that by~\eqref{CsmD}, there is such an edge $CD$ for each $D\in\mathcal D_{i+1}$. So  K\H onig's  lemma  yields  a sequence $(D_i)_{i\in\N}$ with $D_i\in\mathcal D_i$ and $D_i\subseteq D_{i-1}-\partial_v D_{i-1}$ for $i\geq 1$.
 
It is easy to construct a ray $R$ that passes exactly once through each $\partial_v D_i$, and hence there is  an end $\omega\in\bigcap_{i\in\N}\overline D_i$. We claim that
for all $i\in\N$,
\begin{equation}\label{minsepp2}
\text{$\partial_vD_{i+1}$ is an $\subseteq$-minimal  $\partial_vD_i$--$\Omega^G (D_{i+1})$ separator.}
\end{equation}
Then,  $(D_i)_{i\in\N}\rightarrow\omega$. So,  by definition of $\mathcal D_i$, we find that $\df(\omega)<k$, a contradiction to our assumption, as desired.

It remains to show~\eqref{minsepp2}. Let $i\in\N$, and
 observe that by definition of the $S_i$, we know that $\partial_vD_i$ separates the rest of $S_i$ from $D_{i+1}\subseteq D_i-\partial_vD_i$. Hence $N_{S_i}(D_{i+1})\subseteq\partial_vD_i$. As by~\eqref{minsepp},  $\partial_vD_{i+1}$ is an $\subseteq$-minimal $N_{S_i}(D_{i+1})$--$\Omega^G(D_{i+1})$ separator,  and clearly,  $\partial_vD_{i+1}$ is  a  $\partial_vD_i$--$\Omega^G (D_{i+1})$ separator, this implies that $\partial_vD_{i+1}$ is also an $\subseteq$-minimal  $\partial_vD_i$--$\Omega^G (D_{i+1})$ separator, proving~\eqref{minsepp2}.

\end{proof}

We may now use Theorem~\ref{thm:redu} as a black box for translating to infinite locally finite graphs any kind of results from finite graph theory that make assumptions only on the average or minimum degree. For example,
Theorem~\ref{thm:redu} together with Theorem~\ref{thm:minor}/Theorem~\ref{thm:topminor} yields at once
the desired extension of Theorems~\ref{thm:minor} and~\ref{thm:topminor} to locally finite graphs.

\begin{theorem}\label{thm:minorLocFin}
Let $G$ be a locally finite graph. If each vertex and each end of $G$ has (relative) degree at least $ f_1(r)$, then $K^r$ is a minor of~$G$.  If each vertex and each end of $G$ has (relative) degree at least $f_2(r)$, then $K^r$ is a topological minor of~$G$.
\end{theorem}

\medskip

Let us remark that we may not weaken the assumption of Theorem~\ref{thm:minorLocFin} in the following sense.
Denote by $\df '$ the ratio of the edge- and the vertex-degree, that is,  set $\df '(\omega):=\de(\omega)/\dv(\omega)$. 

Now, there is no function $f'$ such that all graphs with $\df '(\omega),d(v)>f'(k)$ for all ends $\omega$ and vertices $v$ contain a complete minor of order $k$. This can be seen by considering the following example (which appeared in a different context in~\cite{hcs}).

\bigskip

\begin{figure}[ht]
      \centering
      \includegraphics[scale=1.1]{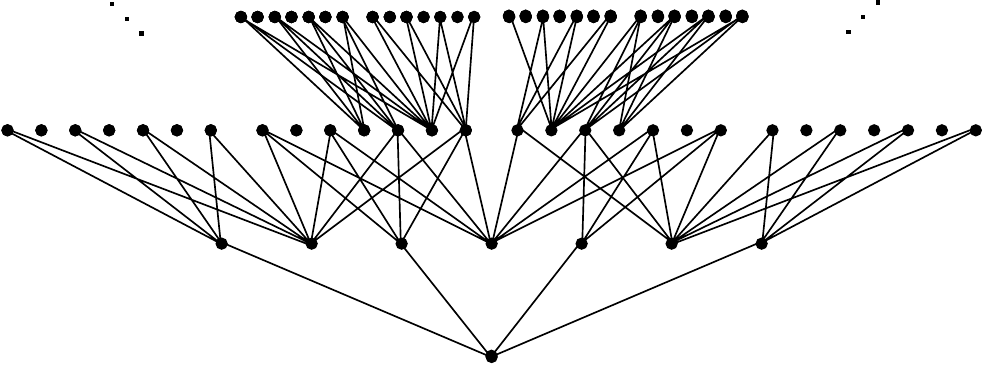}
      \caption{The graph $G_k''$ from Example~\ref{ex:Gk''} for $k=4$.\label{fig:Gk''}}
\end{figure}

\begin{example}\label{ex:Gk''}
Take the infinite tree $T'_k$ with levels $L_0=\{v_0\}, L_1, L_2, \ldots$ where $v_0$ is the root of $T_k$ and each vertex sends $k$ edges to the next level. \\
For each $i\in\N$, consider separately each vertex $x\in L_i$ and its neighbourhood $\{v_1^x,\ldots,v_k^x\}$ in $L_{i+1}$. For $j=1,2,\ldots, k-1$, add a new vertex $w_j^x$ and all edges between $w_j^x$ and $N_{L_{i+2}}(v_j^x)\cup N_{L_{i+2}}(v_{j+1}^x)$ (see Figure~\ref{fig:Gk''}). Call the obtained graph~$G''_k$.\\
Lemma 4.1 of~\cite{hcs} states that $d_e(\omega)\geq k$ and $d_v(\omega)\leq 3$ for all $\omega\in\Omega(G_k'')$. Hence $\df '(\omega)\geq k/3$ for each $\omega\in\Omega(G''_k)$. Clearly, also all vertices of $G''_k$ have degree at least $k$. But, by Lemma 4.3 of~\cite{hcs}, $G''_k$ has no $6$-connected minor, while $k$ may grow  as much as we like.
\end{example}

We finish this section with some problems. First of all, is there a description of the relative degree $\df$ that involves rays instead of sequences of separators? Ideally, this would be similar to the definition of the vertex-/edge-degree.

\begin{problem}
Find an equivalent definition of $\df (\omega)$, in terms of rays of~$\omega$.
\end{problem}

Let us now turn to arbitrary, that is, not necessarily locally finite graphs. First observe that the sequences $(H_i)_{i\in\N}$ that define the relative degree of an end need no longer exist, as for example in $K^{\aleph_0}$. Clearly, the non-existence of these sequences is due to the existence of dominating vertices. Hence, in some way the dominating vertices of an end have to be taken into account for a generalisation of the relative degree notion to arbitrary infinite graphs.

\begin{question}\label{prob:minor}
Is there a natural modification of the relative degree notion that makes an extension of Theorem~\ref{thm:minorLocFin} to arbitrary infinite graphs possible?
\end{question}
 
 A positive answer to this question, at least for graphs with countably many ends, will be given in~\cite{LCM}.

\subsection{Using large girth}\label{sec:girth}

In finite graphs, we can force large complete minors by assuming large girth, and a minimal degree of $3$. More precisely, every finite graph of minimal degree at least $3$ and girth at least $g(k):=8k+3$ has a complete minor of order $k$~\cite{DBook}. 

If we do not take the ends into account, then it is easy to see that this fact does not extend to infinite graphs. Clearly, the $3$-regular infinite tree $T_3$ has infinite girth and no large complete minors, and even if finite girth was required, we might simply add an edge to $T_3$, and still have a counterexample.

But, the ends of our example have end degree $1$ in each of our three end degree notions. Now, we shall see that requiring large minimum vertex- degree at the ends, together with large girth, and minimum degree at least~ $3$ at the vertices, will still not suffice to force large complete minors. 

\begin{figure}[ht]\label{fig:girth}
      \centering
      \includegraphics[scale=0.32]{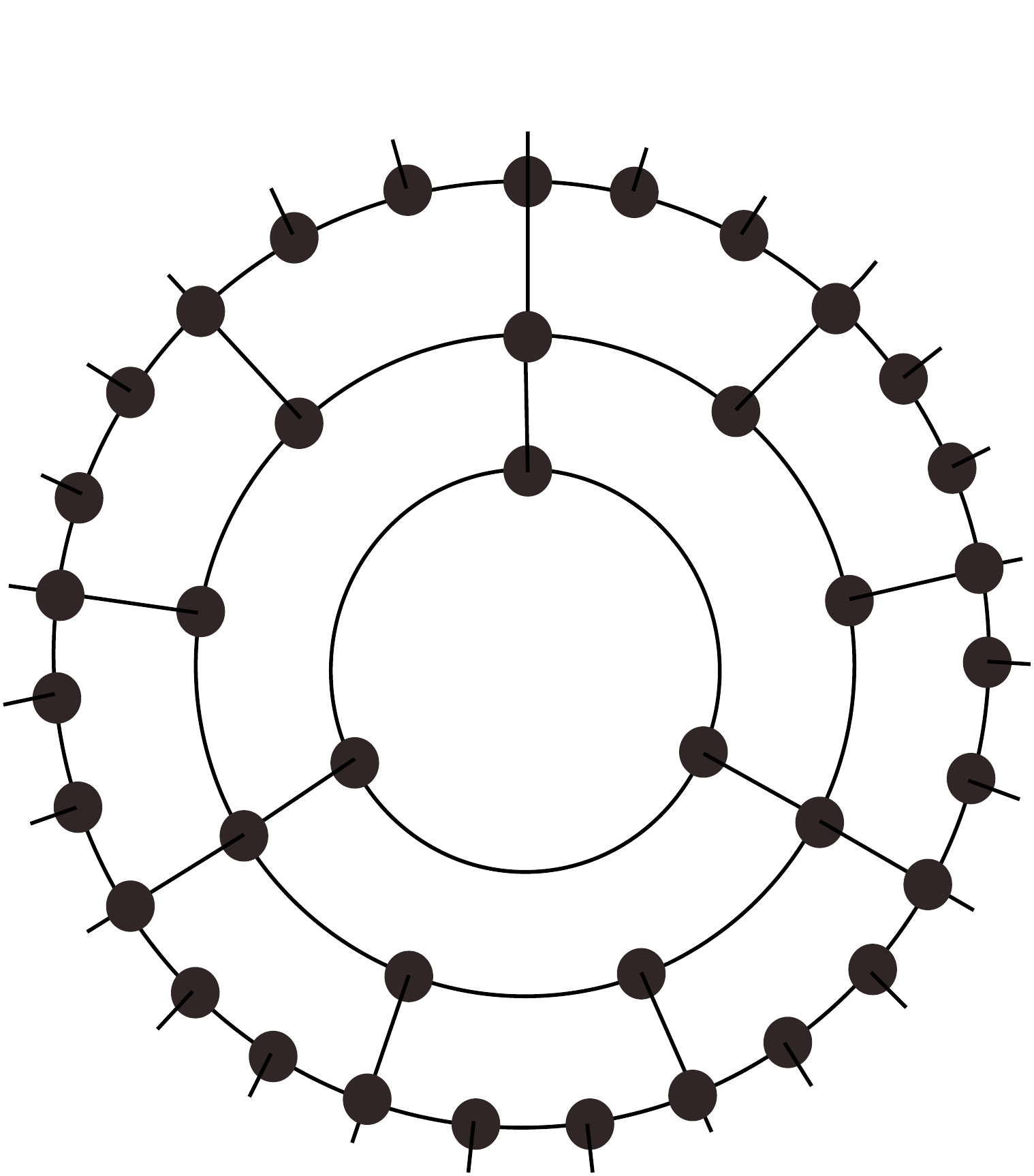}
      \caption{The graph $H_g$ from Example~\ref{ex:girth} for $g=3$.}
\end{figure}

\begin{example}\label{ex:girth}
For all $g\in\N$, we construct a planar graph $H_g$ with finite girth $g$, minimal degree $3$ at the vertices and a unique end, which has infinite vertex-degree.\\
 Take the union of the cycles of length $g^n$, over all $n\in\mathbb N$. We shall add edges between each $C_{g^n}$ and $C_{g^{n+1}}$, one for each vertex in $V(C_{g^n})$, in a way that their new neighbours lie at distance $g$ on $C_{g^{n+1}}$. Clearly, this can be done in a way so that we obtain a planar graph $H_g$ (cf.~Figure~\ref{fig:girth}). Being planar, $H_g$ has no complete minor of order greater than $4$.
\end{example}

However, the relative degree of the end of $H$ is relatively small (in fact, it is~$1$).
Is this a necessary feature of any counterexample?
That is, does every graph of minimum degree $3$ and large girth and without large complete minors have to have an end of small relative degree?
At least the relative degrees cannot be too large:

\begin{proposition}\label{prop:girth}
Every locally finite graph $G$ of minimal degree at least $3$ at the vertices, minimal relative degree at least $r(k)=c_1k\sqrt{\log k}$ at the ends and girth at least $g(k)=8k+3$ has a complete minor of order $k$.
\end{proposition}

\begin{proof}
One may employ the same proof as for finite graphs, as given e.g.~in~\cite{DBook}. The strategy there is to construct first a minor $M$ of $G$ that has large minimal degree, and then apply Theorem~\ref{thm:minor} to $M$. In an infinite graph, we can construct the minor $M$ in exactly the same way, and it is not overly difficult to see that $M$ does not only have large degree at the vertices, but also has at least the same relative degree at the ends as $G$. It suffices to apply Theorem~\ref{thm:minorLocFin} to obtain the desired minor.
\end{proof}

How much can this bound be lowered? May we take $r(k)$ to be constant, even $r(k)=3$? 

\begin{problem}\label{prob:girth}
For $k\in\N$, which is the smallest number $r(k)$ so that every locally finite graph with $d(v)\geq 3$ and $\df (\omega)\geq r(k)$ for all vertices $v$ and ends  $\omega$, and of girth at least $g(k)$ has a complete minor of order $k$?
\end{problem}

%


\bigskip

\section{Minimal $k$-(edge-)connectivity}

\subsection{Edge-minimally $k$-connected graphs}\label{sec:edgeDel}

A $k$-connected graph can be minimal in several ways. The first option that we will investigate here, and which has been most studied until now, is minimality with respect to edge deletion.
Let us call a graph $G$ {\em edge-minimally $k$-connected} if it is $k$-connected but $G-e$ is not, for every $e\in E(G)$.  

Mader~\cite{maderEckenVom} showed that every finite edge-minimally $k$-connected graph $G$ contains at least $|G|/2$ vertices of degree $k$. Halin~\cite{halinUnMin} showed that infinite locally finite edge-minimally $k$-connected  graphs have infinitely many vertices of degree $k$, provided that $k\geq 2$. Mader extended this result to arbitrary infinite graphs.

\begin{theorem}[Mader~\cite{maderUeberMin}]\label{thm:maderEmin}
Let $k\geq 2$ and let $G$ be an infinite edge-minimal\-ly $k$-connected graph. Then the cardinality of the set of those vertices of $G$ that have degree $k$ is $|G|$.
\end{theorem} 

For completeness, let us quickly describe what happens in infinite edge-minimally $1$-connected graphs. Clearly, these are exactly the infinite trees. Thus they do not necessarily have vertices of degree~$1$. But if not, then they must have ends of vertex-degree~$1$. Actually, it is easy to see that they have at least two such points, and unfortunately this is already the best bound for countable trees (because of the double ray). Uncountable trees, however, allow for a version of Theorem~\ref{thm:maderEmin} with ends:

\begin{proposition}
Let $T$ be an  edge-minimally $1$-connected graph (i.e.~a tree) of uncountable order. Then $T$ has $|T|$ vertices of degree~$1$, or $|T|$ ends of vertex-degree~$1$.
\end{proposition}

\begin{proof}
Root $T$ at an arbitrary vertex $r$. Observe that each `leaf', that is,
each vertex/end of (vertex-)degree~$1$ corresponds to a finite or infinite path starting at $r$, and it is easy to see that these paths cover $V(T)$. Hence if $T$ had less than $|T|$ vertices/ends of (vertex-)degree $1$, then $T$ would have order less than $|T|$, a contradiction.
\end{proof}

The proof of Theorem~\ref{thm:maderEmin} relies on the following theorem, which is of interest on its own.

\begin{theorem}[Mader~\cite{maderUeberMin}]\label{thm:fincyc}
Let $k\in\N$, and let $G$ be an edge-minimally $k$-connec\-ted graph. Then each (finite) cycle of  $G$ contains a vertex of degree $k$.
\end{theorem}

In other words, if we delete all vertices of degree $k$ in an edge-minimally $k$-connected graph, we are left with a forest.

It is not overly difficult to see that Theorem~\ref{thm:fincyc}  implies Theorem~\ref{thm:maderEmin} (see~\cite{maderUeberMin}). In order to give an idea, we shall now sketch the easier proof for locally finite $G$. 
The following basic lemma will be useful.
 
 \begin{lemma}\label{lem:raystar}$\!\!${\bf\cite{DBook}}
 Every infinite connected graph contains either a ray or a vertex of infinite degree (or both).
 \end{lemma}

In order to see how the locally finite version of Theorem~\ref{thm:maderEmin} follows from Theorem~\ref{thm:fincyc}, suppose $k$ and $G$ are given as in Theorem~\ref{thm:maderEmin}, and that $G$ is locally finite. If $G$ does not have infinitely many vertices of degree $k$, then by Theorem~\ref{thm:fincyc}, there is a finite non-empty set $S\subseteq V(G)$ so that $F:=G-S$ is a forest. As $k\geq 2$ (by the assumption of Theorem~\ref{thm:maderEmin}), for each $v\in V(F)$ every component of $F-v$ sends at least one edge to $S$. Thus, if $F$ contains a ray $R$, then it is easy to see that there are infinitely many $V(R)$--$S$ edges, contradicting the fact that $G$ is locally finite. Hence, $F$ is rayless, and therefore, by Lemma~\ref{lem:raystar}, has  infinitely many components. These all send edges to $S$, again contradicting the fact that $G$ is locally finite.

\medskip

Mader observed that Theorem~\ref{thm:fincyc} also implies that every subgraph $H$ of a finite edge-minimally $k$-connected graph has vertices of degree at most $k$. In fact, first suppose that $H$ contains a (finite) cycle $C$. Then $C$ is also a cycle in $G$, and thus Theorem~\ref{thm:fincyc} implies that $H$ contains a vertex of degree at most $k$ (in $G$ and thus) in $H$. On the other hand, if $H$ has no finite cycle, then $H$ is a tree and thus has a leaf.

If $G$ and $H$ are infinite then this `leaf' might be an end. Apart from this detail, we may use the same argument for infinite graphs, and thus obtain:

\begin{corollary}\label{prop:Hmin}
 Every subgraph $H$ of an infinite edge-minimally $k$-connec\-ted graph has a vertex of degree at most $k$, or an end of vertex-degree $1$ (in~$H$).
\end{corollary}

\medskip
As mentioned earlier, `infinite cycles', i.e.~circles, play an important role in infinite graph theory.
It is thus natural to ask whether Theorem~\ref{thm:fincyc} extends to circles. It turns out that this is not the case. Infinite circles
 do not necessarily contain vertices of degree~$k$, as can be seen by considering the following example.

\medskip

\begin{figure}[ht]\label{fig:infcyc}
      \centering
      \includegraphics[scale=0.6]{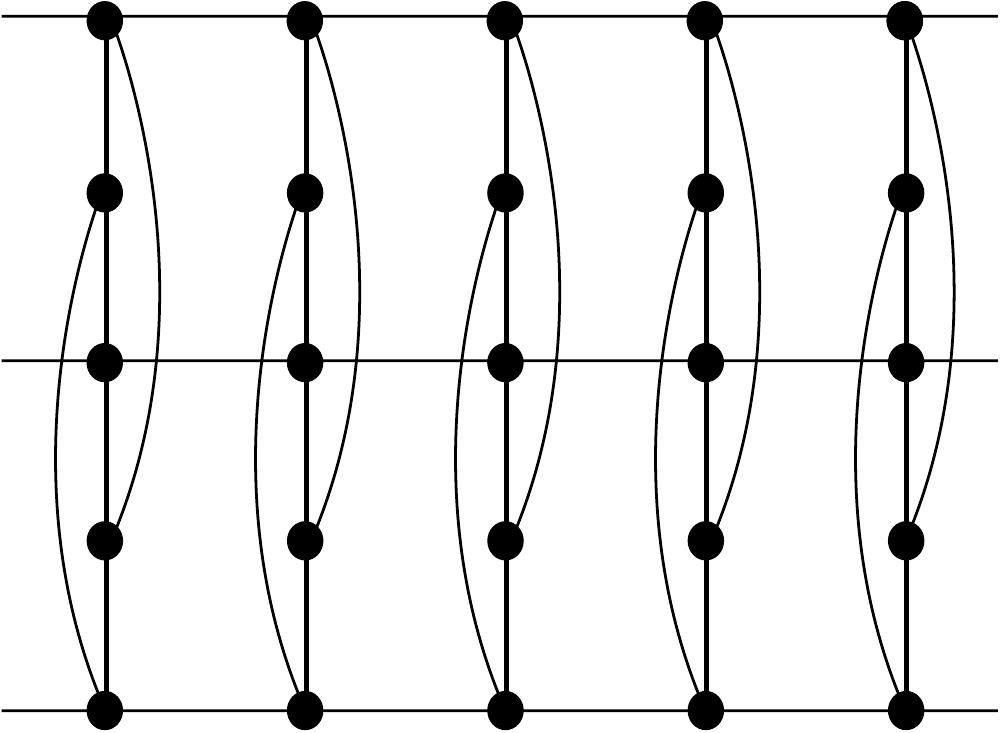}
      \caption{The graph $J_3$ from Example~\ref{ex:infcyc}.}
\end{figure}

\begin{example}\label{ex:infcyc}
Let $k\geq 2$. We define a graph $J_k$ on the vertex set $[2k-1]\times \Z$. Let $J_k$ have the edges $(i,j)(i',j)$ where $i$ mod $2\neq i'$, and the edges $(i,j)(i,j+1)$ for all odd $i$. \\
Clearly, the vertices $(i,j)\in V(J_k)$ have degree $k$ if $i$ is even, and degree $k+1$ if~$i$ is odd. It is easy to see that $J_k$ is $k$-connected. As every edge of $J_k$ is either incident with a vertex of degree $k$, or lies on one of the horizontal $k$-cuts, it follows that $J_k$ is edge-minimally $k$-connected. \\
Now, the vertex set $S:=\{(1,j):j\in \Z\}\cup\{(3,j):j\in\Z\}$ spans a circle in $|J_k|$, while none of the vertices in $S$ has degree $k$.
\end{example}

However, the ends of the graph from Example~\ref{ex:infcyc} have vertex-degree $k$. So, each infinite cycle runs through ends of small degree. This motivates us to ask whether the following infinite version of  Theorem~\ref{thm:fincyc} holds true:

\begin{question}\label{q:cyc}
  Is it true that every (finite or infinite) circle of an infinite edge-minimally $k$-connected graph contains a vertex or an end of (vertex-)\-degree~$k$?
\end{question}

One might be tempted to ask whether something stronger is true, namely, whether {\em all} ends of an edge-minimally $k$-connected graph have vertex-degree $k$.
By Lemma~\ref{lem:converse}, we know that the ends of a locally finite  edge-minimally $k$-connected graph all have vertex-degree at least $k$. So, the question is, can they have larger vertex-degree? Consider the following example to see that the answer is yes.

\begin{figure}[ht]
      \centering
      \includegraphics[scale=0.6]{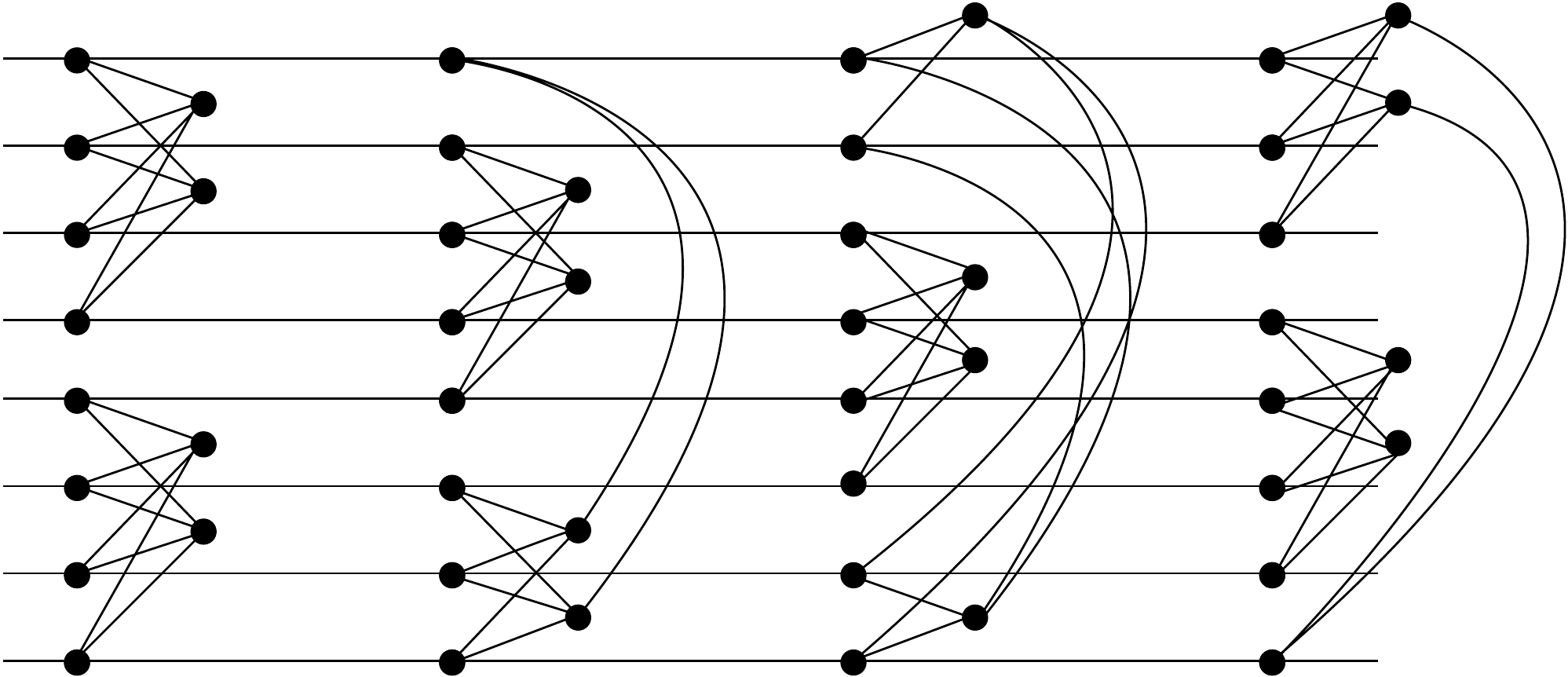}
      \caption{A graph as in Example~\ref{ex:LargeEndDeg} for $k=4$ and $\ell=8$.\label{fig:LargeEndDeg}}
\end{figure}

\begin{example}\label{ex:LargeEndDeg}
We shall construct an edge-minimally $k$-connected locally finite graph whose ends all have vertex-degree $\ell$, for every $k\in\N$ and $\ell\in\N\cup\{ \aleph_0\}$ with $\ell>k$.\\
For each $i\in\Z$, take a copy $H_i$ of $K_{k,k-2}$. Denote by $A_i$ the bigger colour class of $H_i$. Also, take a set $\{R_1,R_2,R_3,\ldots, R_\ell\}$ of double-rays.\\
We now identify the vertices of $A:=\bigcup_{i\in\N}A_i$ with the vertices on the double-rays $R_i$, in a way that each vertex on the $R_i$ gets identified with exactly one vertex of $ A$ and vice versa. We take some care doing this: for each pair of double rays, $R_i=\ldots v^i_{-2}v^i_{-1}v^i_{0}v^i_{1}v^i_{2}\ldots$ and $R_j=\ldots v^j_{-2}v^j_{-1}v^j_{0}v^j_{1}v^j_{2}\ldots$, we can manage that the sets of indices $m$ and $n$ so that $v^i_m$ and $v^j_n$ are mapped to the same $A_r$ are unbounded `in both directions'. More precisely, if $M$ is the set of indices $m$ as in the previous sentence (for some $n$ and some $r$), then $M$ is unbounded from below and from above in $\Z$, and we require the same for set of indices $n$. See Figure~\ref{fig:LargeEndDeg}.
\end{example}

\medskip

We have seen that pretty much is known about infinite edge-minimally $k$-connected graphs. A very important question, however, has not been treated  yet: Does every $k$-connected graph have an edge-minimally $k$-connected subgraph? Or stronger: Does every $k$-connected graph have an edge-minimally $k$-connected spanning subgraph? In finite graphs, the answer is trivially yes: We may simply go on deleting edges as long as the $k$-connectivity is not destroyed. In infinite graphs, this method will not work, and in fact, as Halin~\cite{halinUnMin} pointed out, there are graphs, which have no edge-minimally $k$-connected subgraph at all.

One example of a $2$-connected graph that has no edge-minimally $2$-connected spanning subgraph is the double-ladder $D$, i.e.~the graph on $\{ x_i:i\in\Z\}\cup\{y_i:i\in\Z\}$ with all edges of the form $x_ix_{i+1}$, $y_iy_{i+1}$ or $x_iy_i$. We may delete any  subset $\{x_iy_i:i\in I\}$ of the rungs of $D$ which has the property that $\Z\setminus I$ is unbounded from both above and below, and the graph will stay $2$-connected. But, deleting any other subset of the rungs, our graph will lose its $2$-connectivity. Deleting any other edge but a rung will also destroy the $2$-connectivity. The double-ladder thus has no edge-minimally $2$-connected spanning subgraph. 

\begin{figure}[ht]
      \centering
      \includegraphics[scale=0.6]{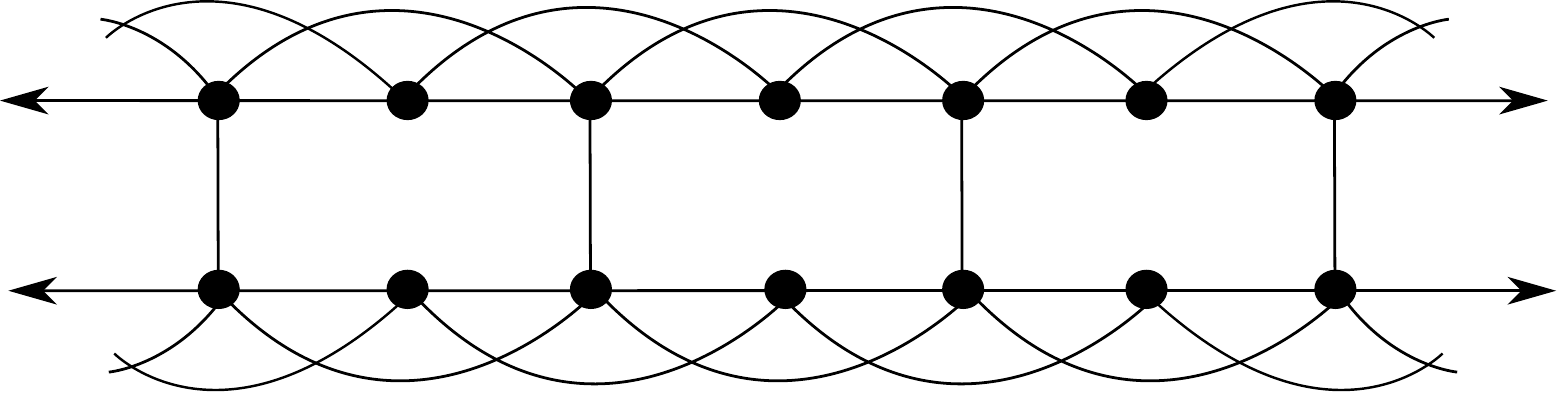}
      \caption{The graph $D_k$. here for $k=2$.\label{fig:D2}}
\end{figure}

Replacing the upper and the lower double-rays of the double-ladder each with the $k$th power of a double-ray, and deleting every second rung, we arrive at a $2k$-connected graph $D_k$ which has no edge-minimally $2k$-connected subgraph at all. Indeed, it is easy to see that any $2k$-connected subgraph of $D_k$ has to contain one and then all vertices of degree $2k$, and therefore all vertices of $D_k$. Now similarly as with the double-ladder, we see that there is no inclusion-maximal subset of the rungs whose deletion leaves the graph $2k$-connected. Hence $D_k$ has no edge-minimally $2k$-connected subgraph.

Coming back to the original example of the double-ladder $D$ (although we could do the same for $D_k$), let us see what happens if we delete all the rungs. Evidently, we arrive at a subgraph $H$ of $D$ that is isomorphic to the (disjoint) union of two double-rays and thus not connected. Viewed as a standard subspace of $|D|$, however, $\overline H$ is still path-connected, moreover, the deletion of any edge, or even of one of the ends, leaves it path-connected. One may thus actually consider the space $\overline H$ to be an edge-minimally $2$-connected standard subspace of $D$. This point of view has been suggested by Diestel~\cite{diestelBanffsurvey}. 

\begin{figure}[ht]
      \centering
      \includegraphics[scale=0.6]{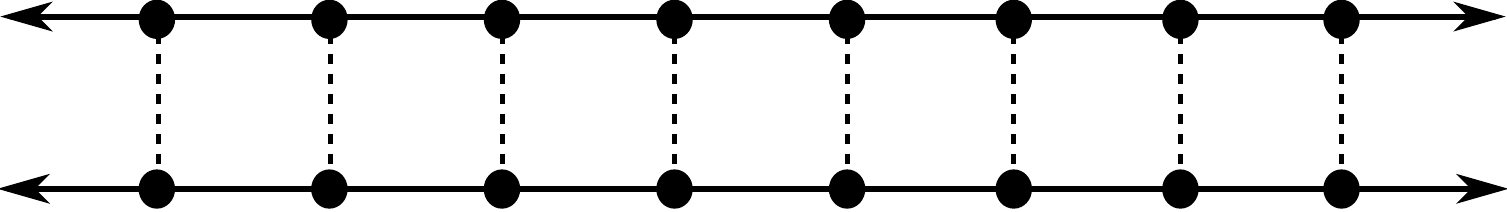}
      \caption{The double-ladder after deleting all its rungs.\label{fig:double-ladder}}
\end{figure}

As we have already discussed in Section~\ref{sec:hcs}, there are two possible notions of $k$-connectivity for standard subspaces: weak and strong $k$-connectivity.
We shall see in Section~\ref{sec:subspace} that an edge-minimally weakly $k$-connected standard subspace exists in every $k$-connected graph, and that the main results on edge-minimally $k$-connected subgraphs carry over to edge-minimally weakly $k$-connected standard subspaces.

\subsection{Vertex-minimally $k$-connected graphs}\label{sec:VDel}\label{sec:vdel}

 Let us now turn to those graphs that are minimally $k$-connected with respect to vertex-deletion. These are the {\em vertex-minimally $k$-connected} graphs, i.e.~those graphs that are $k$-connected but lose this property upon the deletion of any vertex.\footnote{In the literature, these graphs are often called {\em critical} or {\em $k$-critical} graphs, in order to distinguish them from the edge-minimally $k$-connected graphs aka {\em $k$-minimal} graphs. Here,  we chose to speak of edge- and vertex-minimality, in order to make the notation more intuitive.} Clearly, every edge-minimally $k$-connected graph is also vertex-minimally 
 $k$-connected.

Vertex-minimally $k$-connected graphs need no longer contain vertices of degree $k$, but it  has been shown by Chartrand, Lick and Kaugars~\cite{lick} and by Mader~\cite{maderAtome} that finite vertex-minimally $k$-connected graphs necessarily have vertices of degree at most $\frac 32k-1$. Hamidoune~\cite{hamidoune} showed that even two such vertices exist.\footnote{Some authors speak of the bound $\lfloor\frac32k-1\rfloor$. Evidently, this does not make any difference, since the degree is a natural number.}

 The bound $\frac 32 k-1$ on the degree is best possible, as the  following example shows. 
 
\begin{example}\label{ex:lick}
Let $\ell, k\in\N$, and assume that $k$ is even. Take the union of $\ell$ disjoint copies of $K^{k/2}$, which we denote by $H_1,H_2,\ldots,H_\ell$. Add all edges between $H_i$ and $H_{(i+1)mod \ell}$, with $i=1,\ldots \ell$. \\
The obtained graph $O_{k,\ell}$ clearly is  vertex-minimally $k$-connected, and all its vertices have degree $\frac 32 k-1$.
\end{example}

Instead of connecting the first and the last copy of $K^{k/2}$, we may continue infinitely in both directions, thus keeping the minimum degree of the vertices and reducing the connectivity. Replacing the underlying double-ray structure with a tree structure, we see that the degrees of the vertices need not even depend on $k$:

\begin{example}\label{ex:lickInf}
Let $k,r\in\N$, let $T_r$ be the $r$-regular infinite tree, and  for each $v\in V(T_r)$ let $H_v$ be a copy of $K^k$. Take the union of all $H_v$ and  add all edges between $H_v$ and $H_{w}$, if $vw\in E(T_r)$.  Clearly, the obtained graph $T_r(k)$ is vertex-minimally $ k$-connected, and all vertices of  $T_r(k)$ have degree $(r+1)k-1$.
\end{example}

\bigskip

\begin{figure}[ht]
      \centering
      \includegraphics[scale=0.51]{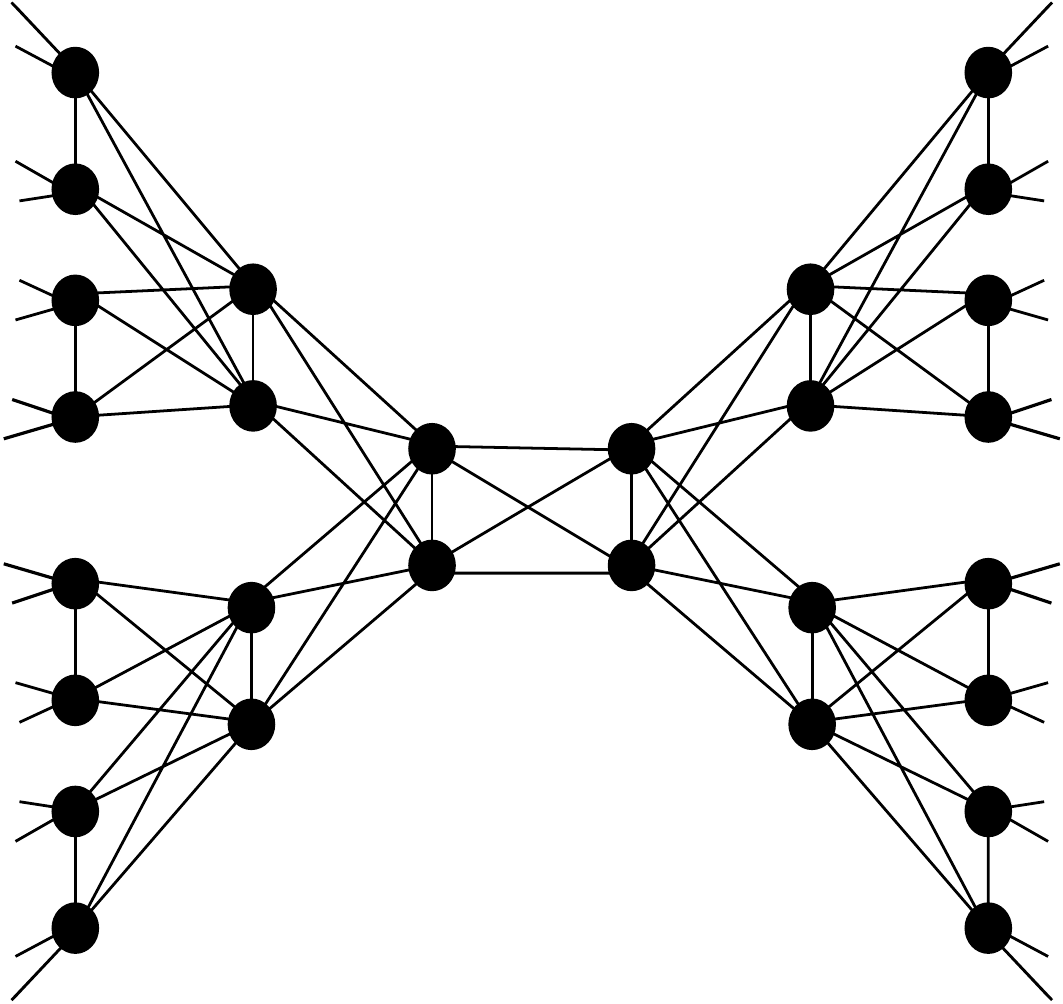}
      \caption{The graph $T_r(k)$ from Example~\ref{ex:lickInf} for $k=2$ and $r=3$.\label{fig:lickInf}}
\end{figure}

However, the vertex-degree of the ends of $T_r(k)$ is $k$. This suggests that an adequate extension of Lick's theorem
to infinite graphs has to allow for ends of small degree. 

And in fact, a first bound is given by Theorem~\ref{thm:hcs}, which implies that every vertex-minimally $k$-connected graph $G$ has a vertex of degree at most $2k(k+1)$ or an end of vertex-degree less than $2k(k+3)$.
But one can do  better:

\begin{theorem}\label{thm:lickInf}$\!\!${\bf\cite[{\rm Theorem 3 (b)}]{minim}}
Let $k\in\N$, and let $G$ be a vertex-minimally $k$-connected graph. Then $G$ has a vertex of degree at most $\frac 32 k-1$ or an end of vertex-degree $\leq k$.
\end{theorem}

One can  improve Theorem~\ref{thm:lickInf} in the spirit of Hamidoune's result mentioned above:

\begin{theorem}\label{thm:lickInf2}$\!\!${\bf\cite[{\rm Theorem 4 (b)}]{minim}}
Let $k\in\N$, and let $G$ be a vertex-minimally $k$-connected graph. Then $|\{\omega\in\Omega(G):d_v(\omega)\leq k\}\cup\{v\in V(G):d(v)\leq \frac 32 k-1\}|\geq 2$.
\end{theorem}

Observe that the bound given by Theorem~\ref{thm:lickInf} is best possible. Indeed, by Lemma~\ref{lem:converse}, the vertex-degree of the ends has to be at least $k$ in a $k$-connected locally finite graph. On the other hand, even if we allow a larger vertex-degree of the ends, we cannot expect a lower  bound on the degrees of the vertices. This is illustrated by the following example.

 \begin{example}\label{ex:Vdel}
  For $k\in 2\N$ and $\ell\in\N\cup\{ \aleph_0\}$, $\ell\geq k$, we construct a vertex-minimally $k$-connected graph $\hat H_{\ell, k}$ whose vertices have degree $\frac 32 k-1$, and whose ends have vertex-degree $\ell$.\\
 Take the disjoint union of $\ell$ double-rays $R_1,\ldots, R_\ell$. For simplicity, assume that $k$ divides $\ell$. For each $i\in\Z$, take $\ell/k$ copies of the graph $O_{k,4}$ from Example~\ref{ex:lick}, and identify the vertices that belong to the first or the last copy of $K^k$ in $O_{k,4}$ with the $i$th vertices the $R_j$. This can be done in a way that the obtained graph has two ends of vertex-degree $\ell$, while the vertices have degree either $3k/2-1$ or $3k/2+1$.
 \end{example}

\bigskip

\begin{figure}[ht]
      \centering
      \includegraphics[scale=0.55]{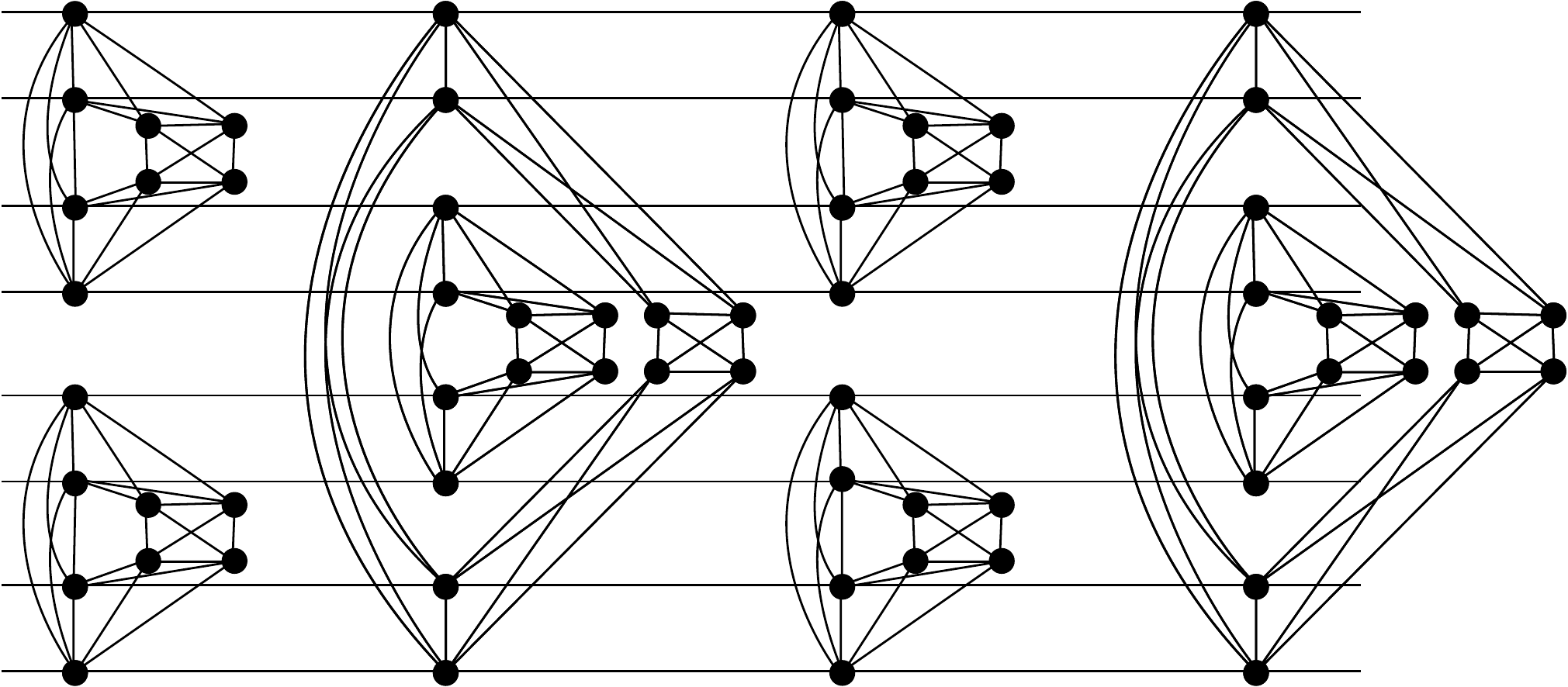}
      \caption{The graph  from Example~\ref{ex:Vdel} for $k=4$ and $\ell=8$.\label{fig:Vdel}}
\end{figure}

Observe that Example~\ref{ex:Vdel} also illustrates the fact that the ends of a vertex-minimally $k$-connected subgraph may all have large vertex-degree (i.e.~independent of $k$), in analogy to Example~\ref{ex:LargeEndDeg}.

\medskip

Finally, we shall ask the same fundamental question as we did for edge-minimally $k$-connected graphs: 

\begin{problem}\label{prob:VMinSubgr}
 Does every $k$-connected graph have a vertex-minimally $k$-con\-nected subgraph?
\end{problem}

In analogy to the edge-minimal case, it is clear that in a finite graph, greedily deleting vertices while not destroying the $k$-connectivity will lead to the desired vertex-minimally $k$-connected subgraph. This need no longer work in infinite graphs, as the following example shows. Take the complete bipartite graph $K_{k,\aleph_0}$. Successively we may delete all vertices of the infinite partition class, at each step maintaining the $k$-connectivity. But after infinitely many steps this procedure will produce a disconnected graph.

`Decontracting' the vertices of the finite partition class of $K_{k,\aleph_0}$ to suitable $k$-connected graphs, we may transform our example to a locally finite one. Indeed, for simplicity assume that $k$ is even, and consider the following example.

\bigskip

\begin{figure}[ht]
      \centering
      \includegraphics[scale=0.51]{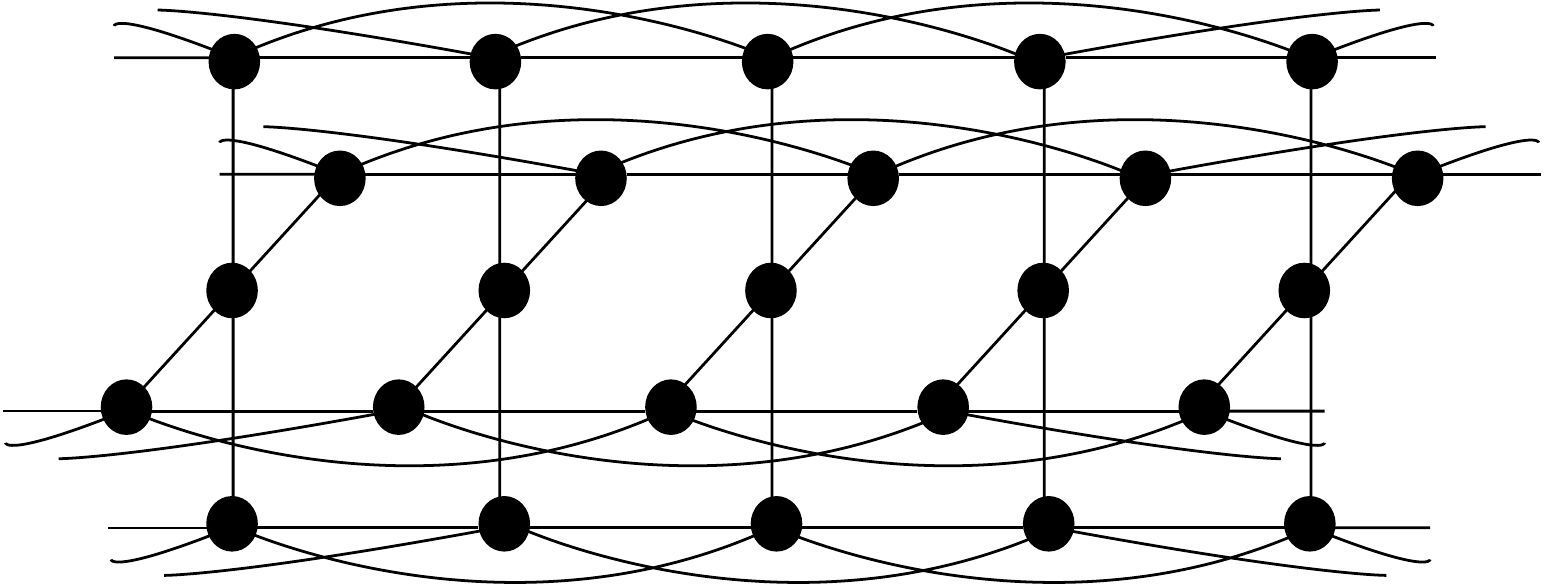}
      \caption{A locally finite $4$-connected graph where greedy deletion of  vertices may fail.\label{fig:greedyVdel}}
\end{figure}

\begin{example}\label{ex:greedyVdel}
 Take $k$ copies of the $\lceil k/2\rceil$th power of the double-ray $R=\ldots v_{-2}v_{-1}v_0v_1v_2\ldots$, and for each $i\in\Z$, let $x_i$ be a new vertex which we shall connect to each of the copies of $v_i$ (see Figure~\ref{fig:greedyVdel}). \\
 It is not difficult to see that this graph is $k$-connected.
\end{example}

Successively deleting all the $x_i$ will at each step will main\-tain the $k$-connec\-ti\-vi\-ty of the graph from Example~\ref{ex:greedyVdel}. But as above, after infinitely many steps we arrive at a disconnected graph.

Note that both our examples contain vertex-minimally $k$-connected subgraphs. In the first example, it is easy to see that the only option for such a subgraph is $K_{k,k}$. 

Our second example, however, has no finite $k$-connected subgraph. This is so because every finite subgraph has  `a last level', most of whose vertices then have degree $<k$ in the subgraph. But nevertheless our graph does have vertex-minimally $k$-connected subgraphs: one such may be obtained by deleting every $(k+1)$st level if $k$ is even, or every $(k+2)$nd level if $k$ is odd.


\subsection{(Induced-)subgraph-minimal $k$-connected graphs}\label{sec:otherMin}

 {\em Subgraph-minimally $k$-connected} graphs, that is, those $k$-connected graphs none of whose proper subgraphs is $k$-connected, might incorporate in some way the concept of minimality of $k$-connected graphs better than  edge- or vertex-minimality do.
Consider the following simple example.
 Take three paths of length at least $3$ and identify their starting vertices, and also identify their endvertices. The obtained graph is both edge- and vertex-minimally $2$-connected. However, our graph has proper $2$-connected subgraphs, namely its cycles.

Instead of subgraph-minimality, we might also consider induced-subgraph-minimality. Define   {\em induced-subgraph-minimally $k$-connec\-ted} graphs as those $k$-connected graphs which have no proper induced $k$-connected subgraph. This is a weaker notion as each subgraph-minimally $k$-connected graph is also induced-subgraph-minimally $k$-connected, and the converse is not true (just consider a long enough square of a cycle to which we add a chord of the cycle).

Clearly, induced-subgraph-minimality implies vertex-minimality. Also, sub\-graph-minimality implies edge-minimality.
All other possible implications do not hold.\footnote{Except for $k=2$, where induced-subgraph-minimality implies edge-minimality.} Finally, edge-minimali\-ty together with induced-subgraph-minimality implies subgraph-minimality. 

So, (that is, since subgraph-minimally $k$-connected graphs are edge-mi\-ni\-mal\-ly $k$-connec\-ted and that induced-subgraph-minimally $k$-connected graphs are vertex-minimally $k$-connected), 
all results of Section~\ref{sec:edgeDel} and Section~\ref{sec:vdel}, respectively, remain true for subgraph-minimally and induced-subgraph-minimally $k$-connected graphs, respectively. In particular, the former have vertices of degree~$k$ on every cycle, and the latter always have vertices of degree at most $\frac 32 k-1$ or ends of vertex-degree at most~$k$. We cannot expect more than this, that is,  we cannot expect to find vertices of lower degree in induced-subgraph-minimally $k$-connected graphs. This is illustrated by Example~\ref{ex:lick}.

\medskip

Now,  although (induced-)subgraph-minimally $k$-connected subgraphs trivially exist in every finite $k$-connected graph, this is no longer clear for infinite graphs. 
As in Example~\ref{ex:greedyVdel}, greedy deletion of vertices might not lead to the desired result. (Nor does greedy edge-deletion, in the case that we aim at subgraph-minimality.)
In fact, as the example from Figure~\ref{fig:D2} shows, infinite $k$-connected graphs need not have subgraph-minimal $k$-connected subgraphs. However, the graphs $D_k$ from Figure~\ref{fig:D2} are themselves induced-subgraph-mi\-ni\-mally $k$-connected.  We ask:

\begin{question}
 Does every $k$-connected graph have an induced-subgraph-mi\-ni\-mally $k$-connected subgraph?
 \end{question}

If not, the following might still be true:

\begin{question}
 Does every $k$-connected graph have a subgraph $H$ such that every induced $k$-connected subgraph of $H$ is isomorphic to $H$?
\end{question}

\medskip

Instead of just allowing the deletion of any kind of subgraphs, one may also consider minimality with respect to deleting certain kinds of subgraphs. 
In this spirit,
Fujita and Kawarabayashi~\cite{FujiKawa} showed that every finite graph that is minimally $k$-connected under the deletion of the endvertices of any edge has a vertex of degree at most $\frac 32 k +1$.  

Moreover, it is conjectured in~\cite{FujiKawa} and proved by Mader~\cite{maderPaths} that there is a function $f:\N\to\N$ such that every finite graph that is minimally $k$-connected under the deletion of any connected subgraph\footnote{Actually, Mader~\cite{maderPaths} proves a stronger result: his graphs are minimally $k$-connected under the deletion of any path of order $\ell$.} of order $\ell$ has a vertex of degree at most $\frac 32 k+f(\ell)$. This extends
the theorem due to Lick et al discussed above in a different direction, and we may ask for the same extensions in infinite graphs:

\begin{question}
Let $G$ be a $k$-connected graph such that for every $xy\in E(G)$, the graph $G-\{x,y\}$ is not $k$-connected. Does $G$ necessarily have a vertex  or an end of (vertex-)degree at most $\frac 32 k+1$?
\end{question}

\begin{question}
Is there a function $f:\N\to\N$ such that the following holds: If $G$ is a $k$-connected graph so that for each connected subgraph $H\subseteq G$ of order $\ell$ the graph $G-H$ is not $k$-connected, then $G$ has a vertex or an end of (vertex-)degree at most $\frac 32 k+f(\ell)$?
\end{question}

We remark that
in finite graphs, the related notion of $(n,k)$-critical graphs has been studied. An $(n,k)$-graph is an $n$-connected graph $G$ that stays $(n-|U|)$-connected upon deletion of any set $U\subseteq V(G)$ of order at most $k$. It has been shown~\cite{maderEndlichkeits} that all $(n,2)$-critical graphs are finite (and they have been determined~\cite{kriesell}). 
Similar holds for $(n,k)$-con-critical graphs, which are finite for all $k>3$ (see~\cite{maderKCon}).

\subsection{Edge-minimally $k$-edge-connected graphs}\label{sec:EDel}

Let us now pose the questions from the previous sections for vertex-/edge-minimally $k$-{\em edge}-connected graphs. In this context, it seems natural to allow for multigraphs instead of (simple) graphs, but then not all results from the finite theory extend, as we shall see below.

 We dedicate this section to {\em edge}-minimally $k$-edge-connected graphs. Finite such graphs have been studied by Lick~\cite{lickline}, who proved that every finite  edge-minimally $k$-edge-connected graph has a vertex of degree $k$. Mader~\cite{maderEdge} proved that unless $k=1$ or $k=3$, there is a constant $c_k$ such that every edge-minimally $k$-edge-connected graph has $c_k|G|$ vertices of degree $k$. For $k=1$ and $k=3$, these constants do not exist: then one can only guarantee for $2$ respectively $4$ vertices of degree $k$.
 
These vertices of small degree need no longer exist in infinite edge-minimally $k$-edge-connected graphs. This can be seen for $k=3$ by considering the square $R^2$ of the double-ray. The graph $R^2$ is edge-minimally $3$-edge-connected, but all its vertices have degree $4$. However, the ends of $R^2$ have edge-degree $3$ (and vertex degree $2$).

\medskip

\begin{figure}[ht]
      \centering
      \includegraphics[scale=.55]{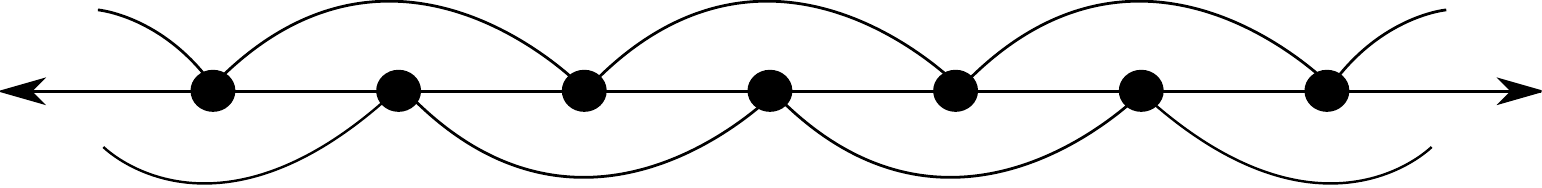}
      \caption{The square of the double-ray.\label{fig:R2}}
\end{figure}

For arbitrary values $k\in\N$, we construct a counterexample as follows.

\begin{example}\label{ex:edgemin}
For $k\in\N$, we construct an edge-minimally $k$-edge-connected graph which has no vertices of degree $k$.\\
 Choose $r\in\N$ and take the $rk$-regular tree $T_{rk}$. For each vertex $v$ in $T_{rk}$, insert edges between the neigbourhood $N_v$ of $v$ in the next level so that $N_v$ spans $r$ disjoint copies of $K^k$. This gives an edge-minimally $k$-edge-connected graph, as one easily verifies. However, the vertices of this graph all have degree at least~$rk$.
\end{example}

But again, in Example~\ref{ex:edgemin} the ends all have edge-degree $k$. This
 gives hope that considering the minimal degrees of the ends might make it possible to extend Lick's theorem to infinite graphs. 
And in fact, this is the case:

\begin{theorem}\label{thm:edgeedge}$\!\!${\bf\cite[{\rm Theorem 3 (c)}]{minim}}
Every edge-minimally $k$-edge-connected graph has a vertex of degree $k$ or an end of edge-degree $\leq k$.
\end{theorem}

This result can be improved in two directions. First, we can guarantee that there are at least two points of small degree, and second, the theorem also holds true for multigraphs (with basically the same proof~\cite{minim}). We thus get:

\begin{theorem}\label{thm:edgeedge2}$\!\!${\bf\cite[{\rm Theorem 4 (c)}]{minim}}
Let $G$ be an edge-minimally $k$-edge-connected multigraph. Then $|\{v\in V(G):d(v)=k\}\cup \{\omega\in \Omega(G):d_e(\omega)\leq k\}|\geq 2$.
\end{theorem}

Observe that in our setting here, it seems more natural to consider the {\em edge}-degree of the ends instead of the vertex-degree (as we are dealing with $k$-{\it edge}-connected graphs). It is also stronger than asking for small  vertex-degree, as the latter (by definition) is bounded from above by the edge-degree.

\medskip

How about an extension of Mader's result mentioned above? Recall that his result states that a positive proportion of the vertices of any finite edge-minimally $k$-edge-connected graph have degree $k$ unless $k=1,3$. 

For infinite graphs $G$, this positive proportion should translate to an infinite set $S$ of vertices and ends of small degree/edge-degree. More precisely, one would wish for a set $S$ of cardinality $|V(G)|$ (or even stronger,  $|S|=|V(G)\cup\Omega(G)|$). 

Also in infinite graphs, we have to exclude the two exceptional values from Mader's result discussed above, $k=1$ and $k=3$. For $k=1$, it is clear that paths  in the finite case, and rays in the infinite case, have only two vertices/ends of (edge-)degree~$1$. For $k=3$, the example of $R^2$ given above illustrates that there are edge-minimally $k$-edge-connected graph that have only two ends of edge-degree~$k$.

\begin{question}\label{q:edge}
For $k\neq 1,3$ does every infinite edge-minimally $k$-edge-connected graph $G$ contain infinitely many vertices or  ends of (edge-)degree $k$? Does $G$ have $|V(G)|$ (or even $|V(G)\cup\Omega(G)|$) such vertices or ends? 
\end{question}
%
%
%

 We remark that Mader's result on the number of vertices of small degree does not hold for multigraphs, no matter whether they are finite or not. For this, it suffices to consider the graph we obtain by multiplying the edges of a finite or infinite path by $k$. This operation results in a multigraph which has no more than the two vertices/ends of (edge-)degree $k$ which were promised by Theorem~\ref{thm:edgeedge2}.

\medskip

Finally, observe that in analogy to the vertex-case, an infinite $k$-edge-connec\-ted graph (or multigraph) need not have a an  edge-minimally $k$-edge-connected spanning subgraph. Again, this can bee seen by considering the double-ladder for $k=2$. Hence it might be interesting to investigate  edge-minimally $k$-edge-connected standard subspaces rather than graphs. 
This question will be shortly addressed in Section~\ref{sec:subspace}.

\subsection{Vertex-minimally $k$-edge-connected graphs}\label{sec:Vedge}

Considering vertex-minimally $k$-edge-connected graphs might seem a little less natural at first sight. Note that, as for $k$-connectivity, the notions `edge-minimally $k$-edge-connected' and `vertex-minimally $k$-edge-connected' are independent in the sense that none implies the other.

 Mader~\cite{maderKritischKanten} showed that every finite vertex-minimally $k$-edge-connected graph contains a vertex of degree $k$, in fact, it contains at least two such vertices. Finite vertex-minimally $k$-edge-connected multigraphs, however, may have arbitrarily large degrees. This can be seen by multiplying each of the edges of $K^k$, and then joining two such modified copies with a maximal matching.

What happens in infinite vertex-minimally $k$-edge-connected graphs? Not only multigraphs, but also simple vertex-minimally $k$-edge-connected graphs need not have vertices of degree~$k$. 

This can already be verified in the double-ladder. In fact, the degrees of the vertices can get arbitrarily large which can be seen in Example~\ref{ex:lickInf}, or even easier in the following modification of it. Replace each vertex of the infinite $r$-regular tree $T_r$ with a copy of the complete graph $K^k$ on $k$ vertices, and add a matching between two of these copies whenever the corresponding vertices of $T_r$ were adjacent (i.e.~we take the product of $T_r$ with $K^k$). But both graphs have ends of vertex-degree $k$. 

This is not a coincidence:

\begin{theorem}\label{thm:vedge}$\!\!${\bf\cite[{\rm Theorems 3 (d) and 4 (d)}]{minim}}
 Let $k\in\N$ and let $G$ be an  infinite vertex-minimally $k$-edge-connected graph. Then $G$ has a vertex of degree~$k$, or an end of vertex-degree at most $k$.\\
Moreover,  $|\{v\in V(G):d(v)=k\}\cup\{\omega\in\Omega(G):d_v(\omega)\leq k\}|\geq 2$.
\end{theorem}

\subsection{Edge-minimally $k$-connected subspaces}\label{sec:edgeDelSubspace}\label{sec:subspace}

We have seen at the end of Section~\ref{sec:edgeDel} that an infinite $k$-connected graph need not contain an edge-minimally $k$-connected subgraph (unless $k=1$). As an example we discussed there the infinite double-ladder $D$. Only the deletion of certain subsets of the rungs of $D$ will leave the graph $2$-connected, but in this way, we will never arrive at an edge-minimally $2$-connected graph.
However, viewing the graph $D^-$ that we obtain by deleting all rungs of $D$ not as a graph on its own, but as a subspace of the space $|D|$, we saw that in fact, we should consider $D^-$ to be a minimally $2$-connected subspace of $D$.

Let us make this idea more precise here. As in Section~\ref{sec:hcs}, for a graph $G$, and a natural number $k$, we call a standard subspace $X$ of $|G|$  that contains at least $k+1$ vertices  {\em weakly $k$-connected} (in $|G|$), if $X-S$ is topologically path-connected for every set $S\subseteq V(G)$ of order less than $k$. We call  $X$  {\em strongly $k$-connected}, if $X-S$ is topologically path-connected for every set $S\subseteq V(G)\cup \Omega(G)$ of order less than $k$.\footnote{It might be interesting to consider also a notion of minimal $k$-connectivity for standard subspaces that lies between weak and strong minimality. One could assign each end a certain weight $w(\omega)$, e.g.~half of its vertex-degree in $X$ (for the definition of end degrees in subspaces see~\cite{degree,DBook}), and then call $X$ minimally $k$-connected if  $X-S$ is topologically path-connected for every set $S\subseteq V(G)\cup \Omega(G)$ with $|S\cap V(G)|+|w(S\cap \Omega(G))|\leq k$. For lack of space, here we do not investigate this promising direction further.} Clearly, strong $k$-connectivity implies weak $k$-connectivity, and it is easy to see that if $X=|G|$, then the usual graph-theoretic $k$-connectivity coincides with strong and weak  $k$-connectivity.

Call $X$ {\em edge-minimally} weakly/strongly $k$-connected, if $X$ is weakly/strongly $k$-connected, but $X-\kreis e$ is not, for every edge $e$ of $G$ with $e\subseteq X$. 
Hence, in particular, if we consider the double-ladder $D$ from above, then 
the closure $\overline{D^-}$ of the subgraph $D^-$ is edge-minimally  strongly (and thus also weakly) $2$-connected in $|D|$.

\medskip

As our motivation for the introduction of these notions was the above-mentioned possible inexistence of edge-minimally $k$-connected subgraphs, the most important question  now is whether every $k$-connected graph $G$ has an edge-minimal weakly or even strongly $k$-connected  standard subspace. If we ask for weakly $k$-connected standard subspaces of locally finite graphs, then the answer is yes. 

\begin{lemma}\label{lem:weak}$\!\!${\bf\cite[{\rm Lemma 3.1}]{diestelBanffsurvey}}
Let $G$ be a locally finite graph, and let $X$ be a weakly $k$-connected standard subspace of $|G|$. Then $X$ has an edge-minimal weakly $k$-connected standard subspace $X'$ which contains $X\cap V(G)$.
\end{lemma}
%
%

Actually, we may obtain such an $X'$ by greedily deleting (the interior of) $k$-connectivity-preseving edges from $X$ (see~\cite{diestelBanffsurvey}).

Note that in particular, Lemma~\ref{lem:weak} implies that every $k$-connected locally finite graph $G$ has an edge-minimal weakly $k$-connected standard subspace which contains $V(G)$.
On the other hand, we do not know whether every  locally finite $k$-connected graph $G$ has an edge-minimal {\em strongly} $k$-connected standard subspace.
Deleting edges greedily we do not necessarily arrive at a strongly $k$-connected standard subspace. 

\begin{figure}[ht]          \centering
      \includegraphics[angle=90,scale=.48]{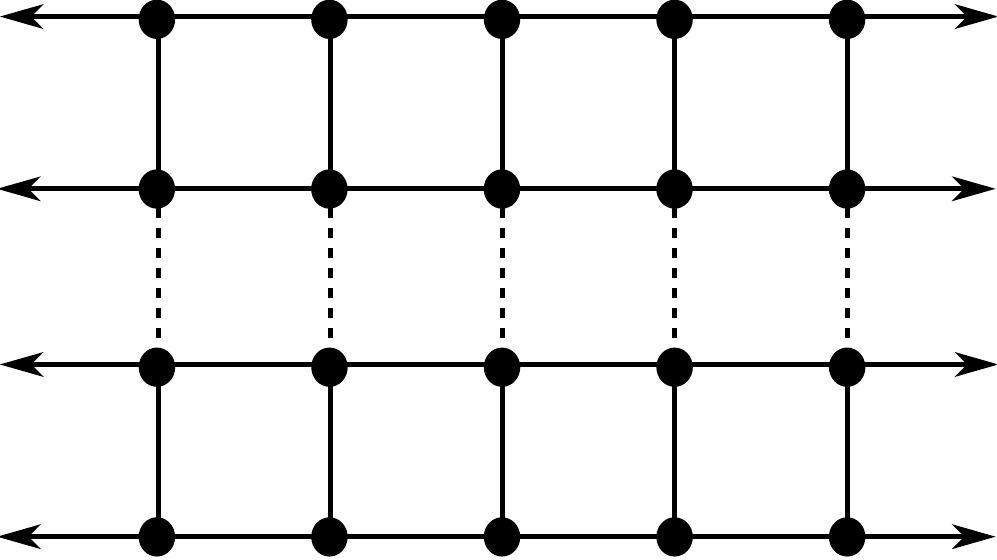}\hspace{1cm}
       \includegraphics[angle=90,scale=.48]{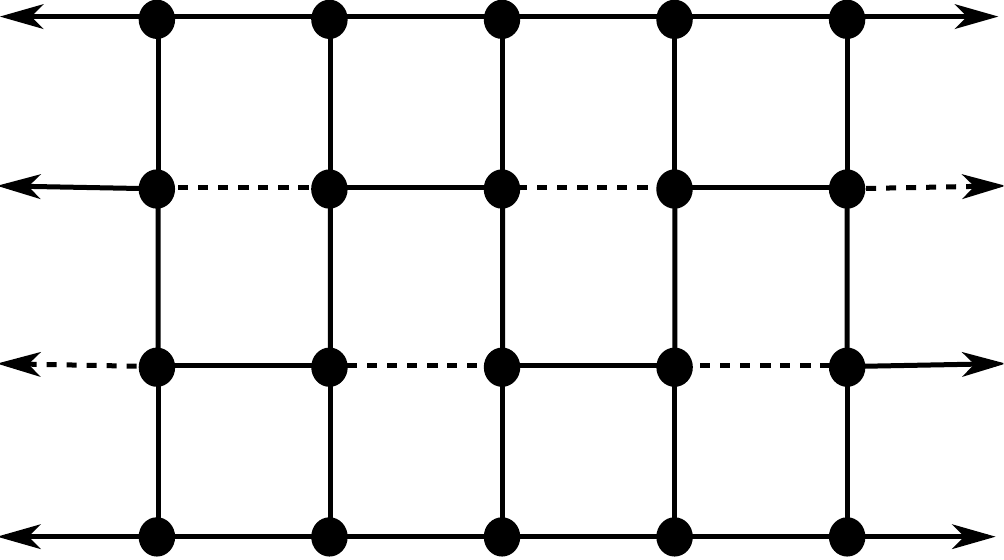}\hspace{1cm}
       \includegraphics[angle=90,scale=.48]{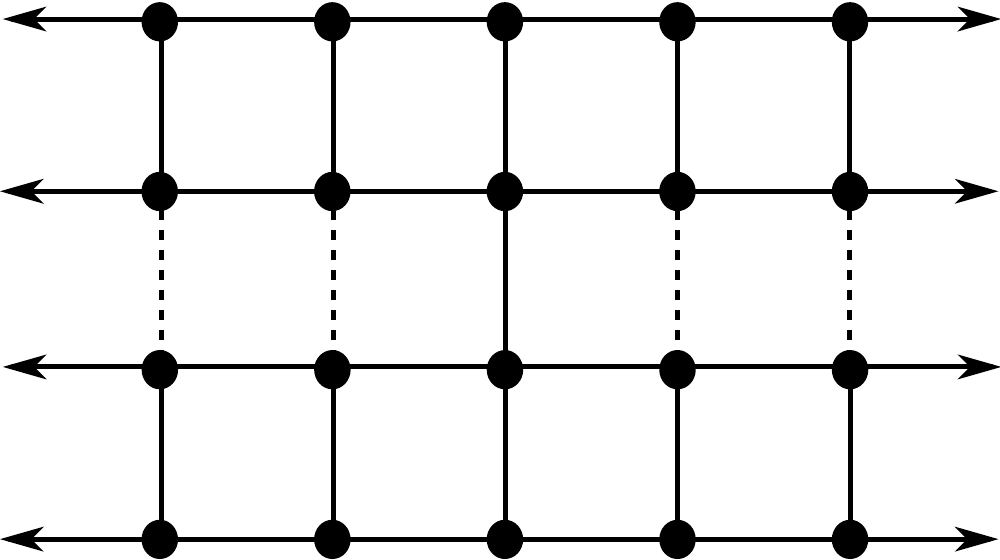}
      \caption{Greedily deleting $k$-connectivity-preserving edges of the $[4]\times \Z$-grid we may arrive at one of these graphs.\label{fig:3zsh}}
\end{figure}

For example, in the $[4]\times \Z$--grid, greedy deletion of $3$-connectivity-preserving edges may lead to many different graphs, three of which we depict in Figure~\ref{fig:3zsh}. The latter two span edge-minimal strongly $3$-connec\-ted standard subspaces of  $[4]\times \Z$, but the first one does not (it is only strongly $2$-connected). Is it  always possible to delete `the right edges'?

\begin{problem}
Let $G$ be a locally finite graph, and let $X$ be a strongly  $k$-connected standard subspace of $|G|$. Is there an edge-minimal strongly $k$-connec\-ted standard subspace $X'\subseteq X$? If so, can $X'$ be chosen so that it contains all of $V(G)\cap X$? 
\end{problem}

Another question would be whether Lemma~\ref{lem:weak} holds for arbitrary (not necessarily locally finite) graphs.

Let us remark that everything said until now in this section remains true if we replace $k$-connected by $k$-edge-connected. In fact, we can define weak and strong $k$-edge-connectivity in the same way as weak and strong $k$-connectivity, and prove an analogue of Lemma~\ref{lem:weak}. Also the problems with strong (edge-) connectivity remain the same.

\medskip

We now show which results from Section~\ref{sec:edgeDel} stay true for standard subspaces. First we shall see that Theorem~\ref{thm:fincyc} carries over.

Recall that we defined the degree $d_X(v)$ of a vertex $v$ in a standard subspace $X$ of the space associated to some graph $G$ as the number of edges $e\in E(G)$ with $e\subseteq X$ that are incident with $v$. 

Then, if $k\geq 2$, every edge-minimally weakly $k$-connected standard subspace of some graph  which contains at least one finite cycle has a vertex of degree $k$ in $X$. We actually have the following stronger result:
 
 \begin{theorem}\label{thm:kappa}
Let $k\geq 2$, let $G$ be a graph, and let $X$ be an edge-minimally weakly $k$-connected standard subspace of $G$ which contains $\kappa$ disjoint finite cycles. Then the cardinality of the set of all vertices of $X$ that have degree $k$ in $X$ is at least $\kappa$. 
\end{theorem}

Theorem~\ref{thm:kappa} follows immediately from  a subspace-version of Theorem~\ref{thm:fincyc}.

\begin{theorem}\label{thm:fincyc2}
Let $k\geq 2$, let $G$ be a graph, and let $X$ be a  weakly $k$-connected standard subspace of $|G|$. Let $C$ be a finite cycle in $G$ such that  $X-\kreis e$ is not weakly $k$-connected for each edge $e$ in $C$. Then $C$ contains a vertex of degree $k$ in $X$.
\end{theorem}

Before we give the proof of Theorem~\ref{thm:fincyc2}, let us remark a few things.
First of all,
observe that the condition that $X$ has a finite cycle is necessary, even if we require $X$ to be strongly $k$-connected. In order to see  this it suffices to consider the following example. 

\begin{example}\label{ex:strongex}
Let $r>k$ be given.
Let $G$ be the product of the infinite $r$-regular tree $T_r$ with a path of length $k$ (i.e.~on $k+1$ vertices). Let $X$ consist of the $k$ copies of $T_r$ plus the end set of $G$. Lemma~\ref{lem:strongex} below asserts that $X$ is edge-minimally strongly $k$-connected. However, all vertices in $V(G)$ have degree $r$ in $X$.
\end{example}

\begin{lemma}\label{lem:strongex}
 The space $X$ from Example~\ref{ex:strongex} is edge-minimally strongly $k$-connected.
\end{lemma}

\begin{proof}
 Since clearly every edge $e\in E(G)$ with $e\subseteq X$ lies in a $k$-cut of $G$, we only have to show that $X$ is strongly $k$-connected. Suppose otherwise. 
 Then there is a set $S\subseteq V(G)\cup\Omega(G)$ with $|S|< k$ so that $X-S$ is not path-connected. Let $x$ and $y$ lie in different path-connected components of $X-S$. We may suppose that $x,y\in V(G)$. 
 
 Since $|S|<k$ there is at least one copy $T^*_r$ of $T_r$ such that $V(T_r^*)\cap S=\emptyset$. Also, as $r>k$, there are rays $R_x$ and $R_y$ starting at $x$ resp.~$y$ which lie completely in $X-S$. Moreover, we can find $R_x$ and $R_y$ such that also their ends lie in $X-S$. Now, $R_x\cup T_r^*\cup R_y\subseteq X-S$ contains an $x$--$y$ arc, a contradiction as desired. 
\end{proof}

The ends of the example just given have vertex-degree $k$, however.
This leads at once to the following question:

\begin{question}\label{q:subsp}
Does every  edge-minimally weakly $k$-connected standard subspace $X$ of an infinite graph $G$ have a vertex or an end of (vertex-)degree $k$?
\end{question}

\bigskip

Observe that Theorem~\ref{thm:fincyc2} also implies a variant of Proposition~\ref{prop:Hmin} for subspaces. In fact, every standard subspace $Y$ of an edge-minimally weakly $k$-connected standard subspace $X\subseteq |G|$ that contains a cycle of $G$ has a vertex of degree $k$, by Theorem~\ref{thm:fincyc2}. On the other hand, if $Y$ has no (finite) cycles, then it may happen that $Y$ has no vertices of degree $k$, and no ends of vertex-degree less than $k$, as is the case in Example~\ref{ex:strongex}. However, $Y$ might have to have ends of vertex-degree $k$ then, so we might repeat Question~\ref{q:subsp} for standard subspaces of $X$.

As for graphs, we do not know whether Theorem~\ref{thm:fincyc2} extends to circles:

\begin{problem}
Let $G$ be a graph and let $X$ be an edge-minimally weakly $k$-connected standard subspace of $|G|$. Does every circle of $G$ with $C\subseteq X$  have a vertex or an end of (vertex-)degree $k$ in $X$? 
What happens if we replace `weakly $k$-connected' with `strongly $k$-connected'?
\end{problem}

\medskip

We dedicate the rest of this section to the proof of Theorem~\ref{thm:fincyc2} which is strongly inspired by Mader's proof~\cite{maderUeberMin}.

\begin{proof}[Proof of Theorem~\ref{thm:fincyc2}]
Suppose $V(C)=\{a_1,a_2,\ldots a_{\ell}\}$ and $C$ has edges $e_i=a_ia_{i+1}$ for $i=1, 2,\ldots ,\ell$ (throughout the proof, we shall understand all indices to be modulo~$\ell$).
For contradiction suppose that for each $i=1, 2,\ldots ,\ell$:
\begin{equation}\label{eq:deg>k}
d_X(a_i)\geq k+1.
\end{equation}

By assumption, for each $i=1, 2,\ldots ,\ell$,  there is a  set $S_i$ in $X_i:=X-\kreis e_i$ so that $X_i-S_i$ is not path-connected. 
For $j\in \{i,i+1\}$ let $C_i^j$ denote the path-connected component of $X_i-S_i$ that contains $a_j$. Set $W_i^j:=(V(G)\cup \Omega(G))\cap C_i^j$.

We claim that

\begin{equation}\label{eq:1}
W_{i+1}^{i+2}-A_i=W_i^{i+1}-B_i.
\end{equation}

where $A_i\in W_{i+1}^{i+2}$ and $B_i\in W_i^{i+1}$ are such that

\begin{equation}\label{eq:2}
|A_i|<|B_i|.
\end{equation}

Then, using~\eqref{eq:1} for $i=1,\ldots,\ell$ we get that 
\begin{align*}
W_1^2-\bigcup_{i=1}^\ell A_i = &\  W_{\ell+1}^{\ell+2}-\bigcup_{i=1}^\ell A_i 
= \ W_1^2-\bigcup_{i=1}^\ell B_i.
\end{align*}

By~\eqref{eq:2}, this means that there is a vertex or end $x$ that lies in  more $A_i$'s than $B_i$'s. But this is impossible, because if  $x$ lies in $A_m$ and $A_{m'}$, say, then by~\eqref{eq:1}, $x$  also lies in some $B_{m''}$ with $m<m''<m'$ (recall that we are viewing all indices modulo~$\ell$). We have thus reached the desired contradiction.

It remains to show the existence of the $A_i$ and $B_i$ satisfying~\eqref{eq:1} and~\eqref{eq:2}. For this, consider the sets
\[
D_i:=W_i^{i+1}\cap W_{i+1}^{i+1}\cap V(G)
\ \text{ and }\ 
\tilde D_i:=W_i^i\cap W_{i+1}^{i+2}\cap V(G).
\]

\begin{figure}[ht]
      \centering
      \includegraphics[scale=0.35]{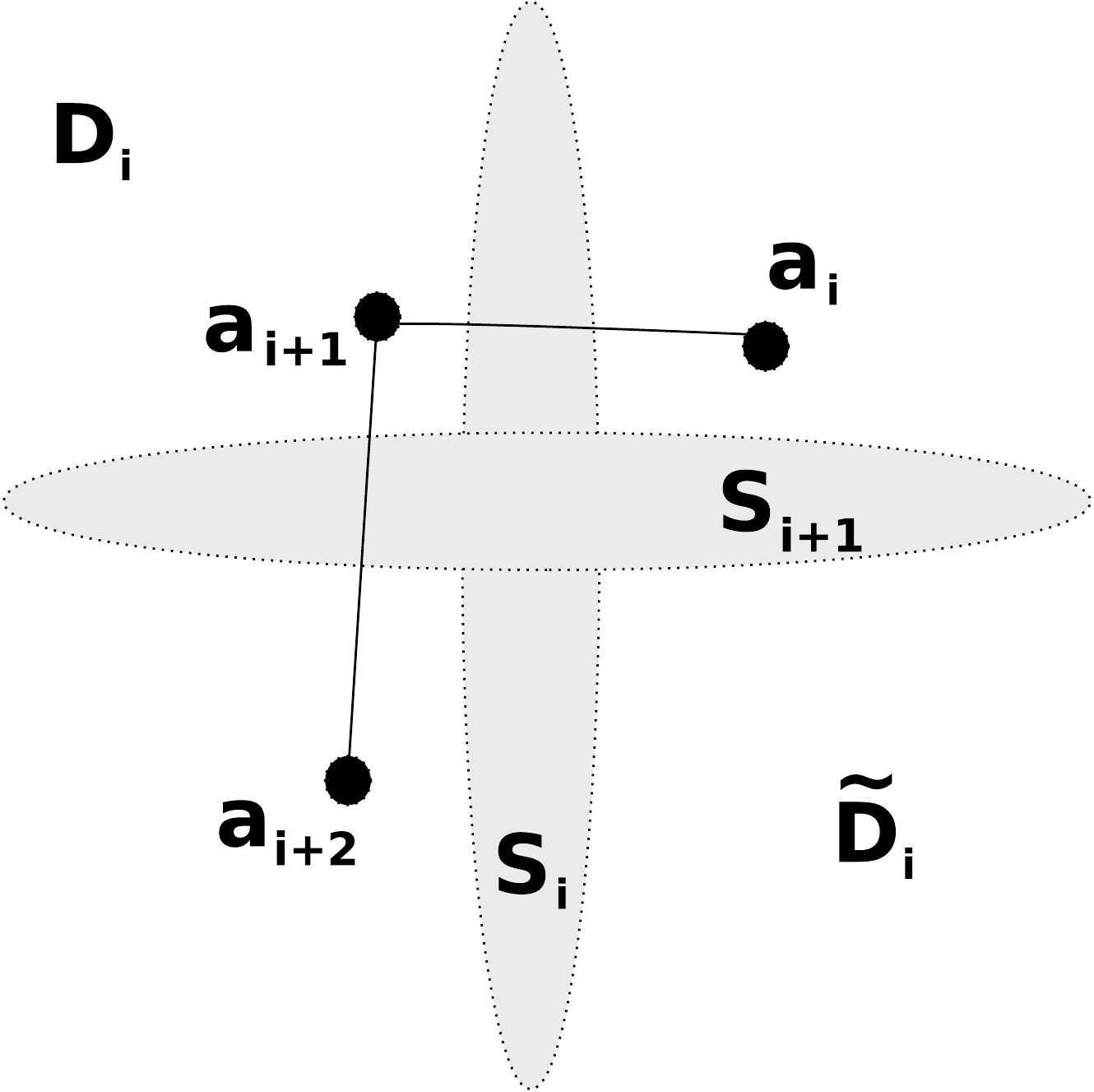}
      \caption{The sets $D_i$ and $\tilde D_i$ in the proof of Proposition~\ref{thm:fincyc2}.\label{fig:Si}}
\end{figure}

For an illustration, see Figure~\ref{fig:Si}. Note that $\tilde D_i$ might be empty. Observe that the neighbourhoods of $D_i$ and $\tilde D_i$ in the subgraph $H\subseteq G$ induced by $X$ satisfy
\[
N_H(D_i)\subseteq (S_{i+1}\cap W_i^{i+1})\cup (S_i\cap S_{i+1}) \cup (S_{i}\cap W_{i+1}^{i+1})\cup \{a_i,a_{i+2}\},
\]
and 
\[
N_H(\tilde D_i)\subseteq (S_{i+1}\cap W_i^{i})\cup (S_i\cap S_{i+1}) \cup (S_{i}\cap W_{i+1}^{i+2}).
\]

Thus,
\begin{equation}\label{eq:kreuz}
|N_H(D_i)|+|N_H(\tilde D_i)|\leq |S_i|+|S_{i+1}|+2 =2k.
\end{equation}

On the other hand, we claim that
\begin{equation}\label{eq:D}
|N_H(D_i)\setminus \{a_i,a_{i+2}\}|\geq k-1.
\end{equation}

Indeed, suppose otherwise. Then $N_H(D_i)\cup \{a_{i+1}\}$ has cardinality at most $k-1$. Hence, as $X$ is weakly $k$-connected, $X-(N_H(D_i)\cup \{a_{i+1}\})$ is path-connected. Since $X-(N_H(D_i)\cup \{\kreis e_i,\kreis e_{i+1}\})$ is not path-connected, this is only possible if
$D_i\setminus\{a_{i+1}\}=\emptyset$. But then $d_X(a_{i+1})<k+1$, a contradiction to~\eqref{eq:deg>k}. This proves~\eqref{eq:D}. 

Combining~\eqref{eq:kreuz} and~\eqref{eq:D} we obtain that 
\begin{align*}
 |N_H(\tilde D_i)|\leq &\ 2k - |N_H(D_i)\setminus \{a_i,a_{i+2}\}|-|\{a_i,a_{i+2}\}|\\ \leq & \ 2k-(k-1)-2\\ =& \ k-1.
\end{align*}
As $X$ is weakly $k$-connected, this implies, similarly as above, that  $\tilde D_i$ is empty.

We set $A_i:=W_{i+1}^{i+2}\cap S_i$, and  $B_i:=W_i^{i+1}\setminus W_{i+1}^{i+2}$. This choice clearly satisfies~\eqref{eq:1}, and for~\eqref{eq:2} it suffices to show that $|A_i|\leq |W_i^{i+1}|\cap S_{i+1}$. So suppose otherwise. Then $|N_H(D_i)\setminus \{a_i,a_{i+2}\}|<|S_i|=k-1$. Hence by~\eqref{eq:deg>k}, $D_i\setminus\{a_{i+1}\}\neq \emptyset$. But then $X-N_H(D_i\setminus\{a_{i+1}\})$ is not path-connected, although $|N_H(D_i\setminus\{a_{i+1}\})|<k$, a contradiction.
This proves~\eqref{eq:2}, and thus completes the proof of the theorem.
\end{proof}

Let us remark that the proof of Theorem~\ref{thm:fincyc2} would also work for edge-minimally weakly $k$-connected spaces $X$ that have the property that $X-\kreis e$ is not {\em strongly} $k$-connected for every edge $e\subseteq X$. In this case, the sets $S_i$ from the proof would consist of vertices {\em and ends}, and instead of the neighbourhoods of $D_i$ and $\tilde D_i$ we would consider certain subsets of $(V(G)\cup\Omega(G))\cap X$. These would  consist of the neighbourhood of $D_i$ and all ends in $\overline D_i$, and the same for~$\tilde D_i$.

\subsection{Other minimally $k$-connected subspaces}\label{otherSubsp}

Let us consider the approach from the previous section for vertex- or subgraph-minimality. That is, we now consider vertex-minimally and (induced-)subgraph-minimally weakly or strongly $k$-connected standard subspaces, defined in the same way as the edge-minimally  weakly or strongly $k$-connected standard subspaces above. Generalising the results and questions of Sections~\ref{sec:vdel} and~\ref{sec:otherMin} we  ask for the existence and the properties of such graphs.

\begin{question}
 Does every $k$-connected graph have a vertex-minimally weakly or even strongly $k$-connected standard subspace? Does every weakly/strongly $k$-connected standard subspace of an infinite graph have a  vertex-minimally weakly/strongly $k$-connected standard subspace?
\end{question}

For locally finite graphs, we can imitate the proof of Lemma~\ref{lem:weak} and obtain a positive answer to the above questions for weakly $k$-connected standard subspaces. Once we have such a space, does it have the vertices of small degree?

\begin{problem}
 Let $X$ be a vertex-minimally weakly/strongly $k$-connected standard subspace of an infinite graph. Does $X$ contain vertices of degree at most $\frac 32 k-1$? 
\end{problem}

For subgraph-minimally weakly $k$-connected standard subspaces, we know from Section~\ref{sec:subspace} that they have vertices of degree $k$ on every (finite) cycle. How about induced-subgraph-minimally weakly/strongly $k$-connected subspaces?

\begin{problem}
 Let $X$ be an induced-subgraph-minimally weakly/strongly $k$-connected standard subspace of an infinite graph. Does $X$ contain vertices or ends of `small' degree? How many?
\end{problem}

If we wish to ask for the existence of these subspaces, the first idea would be to phrase our questions as follows:
 Does every $k$-connected graph have an (induced-)subgraph-minimally weakly/strongly $k$-connected standard subspace? Or, does every $k$-connected standard subspace of an infinite graph have an (induced-) subgraph-minimally weakly or strongly $k$-connected standard subspace?
This, however, might be too strong. Perhaps it would be more natural to ask:

\begin{question}
  Does every weakly $k$-connected standard subspace $X$ of an infinite graph have a weakly $k$-connected standard subspace $Y$ such that all weakly $k$-connected standard subspaces of $Y$ are isomorphic to $Y$? What happens if we replace `weakly $k$-connected' by `strongly $k$-connected'?
\end{question}

\bigskip

\section{Circles, arcs, and forests}\label{sec:3}

\subsection{Hamilton circles}\label{sec:Ham}

Early attempts to generalise results on Hamilton cycles in finite graphs to infinite graphs have been made by Nash-Williams~\cite{NW71}. He considered spanning double-rays as the infinite analogues of Hamilton cycles. This leads to a severe restriction on the class of objects one may study: it is not difficult to see that a graph with a spanning double-ray has at most two ends. So, although some interesting results have been obtained with this notion (Yu~\cite{yuV} proved a conjecture of Nash-Williams that extends Tutte's theorem discussed below to spanning double rays), this  is not quite satisfactory.

From the topological viewpoint on infinite graphs, there is a much more intuitive generalisation of Hamilton cycles. This has been first pointed out by Bruhn. Since circles are the hoemeomorphic image of the unit circle in the space $|G|$ associated to the graph $G$, nothing seems more natural than to adapt this notion and define a {\em Hamilton circle} of a graph $G$ as a circle $C$ in $|G|$ that visits every vertex of $G$. Observe that this definition guarantees that every end of $G$ gets visited `exactly once' by $c$ and that $C$ is a standard subspace of $|G|$.

One of the best known results on Hamilton cycles in finite graphs is due to Tutte. It states that every (finite) $4$-connected  planar graph has a Hamilton cycle. Bruhn (see~\cite{cyclesIntro}) conjectured that this should extend to infinite graphs, using his notion of a Hamilton circle:

\begin{conjecture}[Bruhn]\label{henning}
Every $4$-connected locally finite planar graph has a Hamilton circle.
\end{conjecture}

Partial results on Conjecture~\ref{henning} have been obtained by Bruhn and Yu~\cite{infHamilton} and by Cui, Wang and
Yu~\cite{CWY}.

An extension of Fleischner's theorem on Hamilton cycles to infinite graphs has been conjectured by Diestel~\cite{cyclesIntro} and proved by Georgakopoulos~\cite{HamG}:

\begin{theorem}[Georgakopoulos~\cite{HamG}]\label{agelos}
 Let $G$ be a locally finite $2$-connected graph. Then $G^2$ has a Hamilton circle.
\end{theorem}

It is also shown in~\cite{HamG} that the $3$rd power of any locally finite connected graph has a Hamilton circle. The finite analogue is well-known and not overly difficult to prove. 
%
%
%
%

So, which other assumptions force Hamilton circles in infinite graphs?
Unfortunately, most conditions for Hamilton cycles in finite graphs, like Dirac's theorem, involve degree assumptions that use the order of the graph as a reference. It seems difficult to find a good generalisation of such conditions to infinite graphs. 
 
 Some results from the finite theory, however, use local conditions that do not involve the order of the graph, and thus might allow for extensions to infinite graphs.
Oberly and Sumner~\cite{oberlySumner} showed that every connected locally connected\footnote{A graph is {\em locally connected} if the neighbourhood of each vertex spans a connected subgraph.} claw-free\footnote{A {\em claw-free} graph is one that has no induced subgraph isomorphic to $K_{1,3}$.}
graph of order at least $3$ has a Hamilton cycle. By a result os Asratian~\cite{asratian}, such a graph, if in addition $3$-connected, is even hamilton-connected, which means that every pair of vertices is connected by a Hamiltonian path.
 
We thus feel motivated to ask:

\begin{question}
Does every infinite connected locally connected claw-free graph have a Hamilton circle? 
\end{question}

Defining hamilton-connectivity in the obvious way\footnote{That is, we define a graph $G$ to be {\em hamilton-connected}, if every pair of vertices can be linked by an arc in $|G|$ which contains all of $V(G)$.} for infinite graphs, Bruhn (personal communication) asks the stronger:

\begin{question}
Is every infinite connected locally connected claw-free graph  hamilton-connected?
\end{question}
%

\subsection{Tree-packing and arboricity}\label{sec:tree-p}

A well-known theorem from finite graph theory is the tree-packing
theorem of Nash-Williams~\cite{nashWilliams61} and Tutte~\cite{tutte61}. It states that
a finite multigraph~$G$ has $k$ edge-disjoint spanning trees if and only if every partition of~$V(G)$, into $r\in\mathbb N$ sets say, is crossed by at
  least $k(r-1)$ edges of $G$. (An edge is said to {\em cross} a
given vertex-partition of a graph $G$ if it has its
endvertices in distinct partition sets.)

Disproving a conjecture of
Nash-Williams, Oxley~\cite{Oxley79} constructed a locally finite graph that for $k=2$
satisfies the second condition but has no two edge-disjoint spanning trees.
His graph however, has two edge-disjoint {\em topological spanning trees}, which are defined as topological trees that contain all vertices of the graph.
And in fact, if one replaces the term `spanning tree' from the tree-packing theorem with the term `topological spanning tree', then the theorem does extend to locally finite graphs. This has been shown by Bruhn and Diestel (see~\cite{DBook}), building on work of Tutte. 

\begin{theorem}$\!\!${\bf\cite{DBook}}
\label{thm:locFinTreePack}
For a locally finite multigraph $G$ the following statements are
equivalent:
\begin{enumerate}[\rm (i)]
\item $|G|$ contains $k$ edge-disjoint topological spanning trees;
\item every partition of~$V(G)$, into $r\in\mathbb N$ sets say, is crossed by at
  least $k(r-1)$ edges of $G$.
\end{enumerate}
\end{theorem}

A related result is Nash-Williams' arboricity theorem, which states that a graph  is the
edge-disjoint union of at most $k$ forests, if no
set of $\ell$ vertices induces more than $k(\ell-1)$ edges. A standard compactness argument shows that Nash-Williams' arboricity theorem extends to infnite graphs, if we ask for traditional forests, i.e.~subgraphs of $G$ that have no finite cycles.
But, having taken the topological viewpoint, one should want more. In fact, it is now natural to require that the graph decomposes into  topological forests. This, however, is false without additional constraints.

In fact, for any $k\in\N$, there are examples of graphs which satisfy the condition of local sparseness, but do not decompose into as few topological forests as desired. It suffices to take one copy of $K^{2k}$ for each $n\in\Z$, and identify, for each $n\in\N$, one vertex of the $n$th copy with one vertex of the $(n+1)$th copy, not using any vertex twice. Then add, for each $n\in\N$, an edge between two not yet used vertices of the $n$th and the $(-n)$th copy.
 It is not difficult to see that the obtained graph is the
edge-disjoint union of $k$ ordinary forests, and hence satisfies Nash-Williams' 
condition that no set of $\ell$ vertices spans more than $k(\ell -1)$ edges. But 
any partition of $G$ into $k$ forests induces such a partition in each copy 
of $K^{2k}$, ie.~into spanning trees of $K^{2k}$. Each of these
contains a $v$--$w$ path, so each of our $k$ forests contains a double ray and 
thus an infinite cycle.

So what goes wrong in this counterexample? In fact, our situation is reciprocal to the one in the beginning of the paper, when we tried to get from local density (implied by large vertex degrees) to global density. This would only work if we required denseness at the ends as well.
Analogously, here we have to impose a sparseness condition on the ends of the graph. This sparseness condition can be expressed in terms of the vertex-degree:

\begin{theorem}\label{thm:arbo}$\!\!${\bf\cite{arbo}}
Let $k\in\mathbb N$, and let $G$ be a locally finite graph in which no
set of $\ell$ vertices induces more than $k(\ell-1)$ edges. Further, let
every end of $G$ have vertex-degree $<2k$. Then $|G|$ is the
edge-disjoint union of at most $k$ topological forests in $|G|$.  
\end{theorem}

Although, as we have seen in the example above, the bound of~$2k$ in Theorem~\ref{thm:arbo} cannot be 
reduced, the theorem has no direct converse: a partition into $k$ 
topological forests does not force all end degrees to be small. The 
$\mathbb N\times\mathbb N$ grid, for example, is an edge-disjoint union of two topological forests (its horizontal vs.\ its vertical edges), but its 
unique end has infinite vertex-degree.

It would be interesting to investigate whether Theorem~\ref{thm:locFinTreePack} and~\ref{thm:arbo} extend to subspaces. In the former theorem, we then have to replace the term `crossing edges' with something like `crossing arcs'. This seems to be necessary as can be seen by considering once again the double-ladder $D$ and its subgraph $D^-$ which is obtained by deleting all the rungs. Now a partition of $\overline{D^-}$ into the two sets corresponding to the two double-rays contained in $D^-$ (and putting the ends anywhere) has no crossing edges. However, $\overline{D^-}$ does not contain a topological spanning tree of $|D|$. This motivates 
the following definition.

For a partition of $V(G)\cup\Omega(G)$ of a graph $G$, an arc $A\subseteq |G|$ is said to cross the partition, if it has its endpoints in different partition sets $P_1$ and $P_2$, and furthermore, $A\cap (V(G)\cup\Omega(G))\subseteq P_1\cup P_2$.

\begin{problem}
Let $G$ be a locally finite multigraph $G$, and let $X\subseteq |G|$ be a standard subspace. Are the following statements equivalent?
\begin{enumerate}[\rm (i)]
\item $X$ contains $k$ edge-disjoint topological spanning trees of $G$;
\item every partition of~$V(G)\cup\Omega(G)$, into $r\in\mathbb N$ sets say, is crossed by at
  least $k(r-1)$ edge-disjoint arcs  $A\subseteq X$ .
\end{enumerate}
\end{problem}

For a version of Theorem~\ref{thm:arbo} for subspaces, we use the definition of the vertex-degree in standard subspaces, which can be found at the end of Section~\ref{sec:hcs}.

\begin{problem}
 Let $k\in\mathbb N$, let $G$ be a locally finite graph, and let $X$ be a standard subspace of $|G|$.\\ If no
set of $\ell$ vertices of $G$ induces more than $k(\ell-1)$ edges $e$ with $e\subseteq X$ and furthermore,
every end of $G$ has vertex-degree $<2k$ in $X$, is then $X$ is 
edge-disjoint union of at most $k$ topological forests in $|G|$?  
\end{problem}

\subsection{Connectivity-preserving arcs and circles}\label{sec:con-pres}

There are a few very interesting conjectures about connectivity-preserving paths and cycles in finite graphs. The most famous among these is a conjecture of Lov\'asz (see~\cite{ThomGraphDec}):

\begin{conjecture}[Lov\'asz]\label{conj:lov}
 There is a function $f(k)$ such that for every finite $f(k)$-connected graph $G$, and every pair of vertices $v$, $w$ of $G$ there is an induced $v$--$w$ path $P$ such that $G-V(P)$ is $k$-connected.
\end{conjecture}

The conjecture can equivalently be stated as follows: There is a function $f(k)$ such that for every finite $f(k)$-connected graph $G$ and every edge $e$ of $G$ there is an induced cycle $C$ containing $e$ so that $G-V(C)$ is $k$-connected.
Lov\'asz also conjectured that if we do not insist on prescribing an edge which the cycle has to contain, then $f(k)=k+3$. This has been verified by
Thomassen~\cite{thomassenNonSep}: every finite $(k+3)$-connected graph has an induced cycle $C$ so that  the deletion of $V(C)$ results in a $k$-connected graph. 

A weakening of Conjecture~\ref{conj:lov} has been conjectured by Kriesell and proved by Kawarabayashi et al~\cite{weakLovasz}. It states that there is a function $f(k)$ so that for every $f(k)$-connected graph $G$ and for every edge $e$ of $G$ there is an induced cycle $C$ of $G$ with $e\in E(C)$ such that $G-E(C)$ is $k$-connected.

These results (and Lov\'asz' conjecture, if true) do not carry over to infinite graphs, if we ask for connectivity-preserving cycles that are {\em finite}. The reason is that there are infinite (even locally finite) graphs of arbitrarily large connectivity whose cycles are all separating. More precisely, for every $k\in\N$ there is a $k$-connected locally finite graph such that for each cycle $C$ of $G$ we have that both $G-V(C)$ and $G-E(C)$ are  disconnected. Such graphs have been constructed by Aharoni and Thomassen~\cite{AT}.

Will it help to consider circles instead of finite cycles, and arcs instead of finite paths? Some of the following problems have been posed already in~\cite{cyclesIntro}. Call an arc $A$ or a circle $C$ {\em induced} if $e\subseteq A$ for each edge $e$ with both endpoints in $A$ resp.~$C$.

\begin{problem}
 Is there a function $f(k)$ such that for every infinite $f(k)$-connected graph $G$, and every pair of vertices $v$, $w$ of $G$ there is an induced $v$--$w$ arc $A$ whose deletion leaves a  strongly/weakly $k$-connected subspace of $|G|$? If so, may we also prescribe ends to be the starting points/endpoints of $A$?
\end{problem}

\begin{problem}
 Does every infinite $(k+3)$-connected graph have an induced circle $C$ so that  the deletion of $V(C)$ results in a weakly/strongly $k$-connected subspace of~$|G|$?
\end{problem}

\begin{problem}
 Is there  a function $f(k)$ so that for every edge $e$ of an infinite $f(k)$-connected graph  $G$ there is an induced circle $C$ of $G$ which contains $e$ such that the deletion of the edges of $C$ results in a strongly/weakly $k$-connected subspace of~$|G|$?
\end{problem}

Stronger versions of these problems can be obtained by replacing the graph $G$ with a standard subspace $X$.

An analogue of Thomassen's result for edge-connectivity also holds. Indeed, Mader~\cite{maderKreuz} showed that every finite $(k+2)$-edge-connected graph contains an induced cycle $C$ such that the deletion of $E(C)$ leaves a $k$-edge-connected graph. 

Recently, this has been extended to infinite graphs by Bruhn, Diestel and Pott~\cite{dualTrees}, using the notion of weak $k$-edge-connectivity (which one defines analogously to weak $k$-connectivity).
More might be true:

\begin{problem}
 Is there  a function $f(k)$ so that for every edge $e$ of an $f(k)$-edge-connected graph  $G$ there is an (induced) circle $C$ of $G$ which contains $e$ such that the deletion of the edges of $C$ results in a strongly/weakly $k$-edge-connected subspace of~$|G|$?
\end{problem}

\vskip1.8cm

\bibliographystyle{plain}
\bibliography{graphs}

\begin{thebibliography}{10}

\bibitem{AT}
R.~Aharoni and C.~Thomassen.
\newblock Infinite highly connected digraphs with no two arc-disjoint spanning
  trees.
\newblock {\em J.~Graph Theory}, 13:71--74, 1989.

\bibitem{asratian}
A.S. Asratian.
\newblock Every $3$-connected, locally connected claw-free graph is
  hamilton-connected.
\newblock {\em J.~Graph Theory}, 23:191--201, 1996.

\bibitem{BBExtGT}
B.~Bollob\'as.
\newblock {\em Extremal {G}raph {T}heory}.
\newblock Academic Press London, 1978.

\bibitem{BBModern}
B.~Bollob\'as.
\newblock {\em Modern {G}raph {T}heory}.
\newblock Springer-Verlag, 1998.

\bibitem{dualTrees}
H.~Bruhn, R.~Diestel, and J.~Pott.
\newblock Dual trees can share their ends.
\newblock Preprint 2009.

\bibitem{degree}
H.~Bruhn and M.~Stein.
\newblock On end degrees and infinite circuits in locally finite graphs.
\newblock {\em Combinatorica}, 27:269--291, 2007.

\bibitem{infHamilton}
H.~Bruhn and X.~Yu.
\newblock Hamilton circles in planar locally finite graphs.
\newblock {\em SIAM. J. Discrete Math.}, 22:1381--1392, 2008.

\bibitem{lick}
G.~Chartrand, A.~Kaugars, and D.~Lick.
\newblock Critically $n$-connected graphs.
\newblock {\em Proc.~Am.~Math.~Soc.}, 32:63--68, 1972.

\bibitem{CWY}
Q.~Cui, J.~Wang, and X.~Yu.
\newblock Hamilton circles in infinite planar graphs.
\newblock {\em J. Comb. Theory Ser. B}, 99(1):110--138, 2009.

\bibitem{diestelBanffsurvey}
R.~Diestel.
\newblock Locally finite graphs with ends: a topological approach.
\newblock Preprint 2009 ({H}amburger {B}eitr\"age zur {M}athematik). \\Note: A
  part of this survey will appear in Discr.~Math.'s special issue on infinite
  graphs, another part in C.~Thomassen's birthday volume (also Discr.~Math.).

\bibitem{cyclesIntro}
R.~Diestel.
\newblock The cycle space of an infinite graph.
\newblock {\em Comb.,\ Probab.\ Comput.}, 14:59--79, 2005.

\bibitem{DBook}
R.~Diestel.
\newblock {\em Graph Theory \emph{(4th edition)}}.
\newblock Springer-Verlag, 2010.

\bibitem{cyclesI}
R.~Diestel and D.~K{\"u}hn.
\newblock On infinite cycles {I}.
\newblock {\em Combinatorica}, 24:69--89, 2004.

\bibitem{cyclesII}
R.~Diestel and D.~K{\"u}hn.
\newblock On infinite cycles {II}.
\newblock {\em Combinatorica}, 24:91--116, 2004.

\bibitem{tst}
R.~Diestel and D.~K{\"u}hn.
\newblock Topological paths, cycles and spanning trees in infinite graphs.
\newblock {\em Europ.\ J.\ Combinatorics}, 25:835--862, 2004.

\bibitem{frankHofC}
A.~Frank.
\newblock Connectivity and network flows.
\newblock In {\em Handbook of {C}ombinatorics, Vol.\ 1}, pages 111--177.
  Elsevier, Amsterdam, 1996.

\bibitem{FujiKawa}
S.~Fujita and K.~Kawarabayashi.
\newblock Connectivity keeping edges in graphs with large minimum degree.
\newblock {\em J.~Combin.\ Theory (Series B)}, 98:805--511, 2008.

\bibitem{HamG}
A.~Georgakopoulos.
\newblock Infinite {H}amilton cycles in squares of locally finite graphs.
\newblock {\em Adv.\ Math.}, 220:670--705, 2009.

\bibitem{halin65}
R.~Halin.
\newblock {\"U}ber die {M}aximalzahl fremder unendlicher {W}ege in {G}raphen.
\newblock {\em Math.\ Nachr.}, 30:63--85, 1965.

\bibitem{halinUnMin}
R.~Halin.
\newblock Unendliche minimale $n$-fach zusammenh\"angende {G}raphen.
\newblock {\em Abh.~Math.~Sem.~Univ.~Hamburg}, 36:75--88, 1971.

\bibitem{hamidoune}
Y.~O. Hamidoune.
\newblock On critically $k$-connected graphs.
\newblock {\em Disc.\ Math.}, 32:257--262, 1980.

\bibitem{weakLovasz}
K.~Kawarabayashi, O.~Lee, B.~Reed, and P.~Wollan.
\newblock A weaker version of {L}ov\'{a}sz' path removal conjecture.
\newblock {\em J. Comb. Theory Ser. B}, 98(5):972--979, 2008.

\bibitem{kriesell}
M.~Kriesell.
\newblock Mader's conjecture on extremely critical graphs.
\newblock {\em Combinatorica}, 26:277--314, 2006.

\bibitem{endsBerniElmar}
B.~Kr\"on and E.~Teufl.
\newblock Ends -- {G}roup-theoretical and topological aspects.
\newblock Preprint 2009.

\bibitem{lickline}
D.~R. Lick.
\newblock Minimally $n$-line connected graphs.
\newblock {\em J.~Reine Angew.~Math.}, 252:178--182, 1972.

\bibitem{maderSym}
W.~Mader.
\newblock {\"U}ber den {Z}usammenhang symmetrischer {G}raphen.
\newblock {\em Archiv der Math.}, 21:331--336, 1970.

\bibitem{maderAtome}
W.~Mader.
\newblock Eine {E}igenschaft der {A}tome endlicher {G}raphen.
\newblock {\em Arch.~Math.}, 22:333--336, 1971.

\bibitem{maderMin}
W.~Mader.
\newblock Minimale $n$-fach zusammenh\"angende {G}raphen.
\newblock {\em Math.~Ann.}, 191:21--28, 1971.

\bibitem{maderEckenVom}
W.~Mader.
\newblock Ecken vom {G}rad n in minimalen n-fach zusammenh\"angenden {G}raphen.
\newblock {\em Arch. Math. (Basel)}, 23:219--224, 1972.

\bibitem{maderUeberMin}
W.~Mader.
\newblock {\" U}ber minimale, unendliche n-fach zusammenh\" angende {G}raphen
  und ein {E}xtremalproblem.
\newblock {\em Arch. Math. (Basel)}, 23:553--560, 1972.

\bibitem{maderKreuz}
W.~Mader.
\newblock Kreuzungsfreie a,b-{W}ege in endlichen {G}raphen.
\newblock {\em Abh. Math. Sem. Univ. Hamburg}, 42:187--204, 1976.

\bibitem{maderEndlichkeits}
W~Mader.
\newblock Endlichkeitss\" atze f\"ur $k$-kritische {G}raphen.
\newblock {\em Math.Ann.}, 229(2):143--153, 1977.

\bibitem{maderPaths}
W.~Mader.
\newblock Paths in graphs, reducing the edge-connectivity only by two.
\newblock {\em Graphs and Combinatorics}, 1(1):81--89, 1985.

\bibitem{maderKritischKanten}
W.~Mader.
\newblock Kritisch $n$-fach kantenzusammenh\"angende {G}raphen.
\newblock {\em J.~Combin.\ Theory (Series B)}, 40:152--158, 1986.

\bibitem{maderEdge}
W.~Mader.
\newblock On vertices of degree $n$ in minimally $n$-edge-connected graphs.
\newblock {\em Combinatorics, Probability {\&} Computing}, 4:81--95, 1995.

\bibitem{maderKCon}
W.~Mader.
\newblock On $k$-con-critically $n$-connected graphs.
\newblock {\em J.~Combin.\ Theory (Series B)}, 86:296--314, 2002.

\bibitem{nashWilliams61}
C.St.J.A. Nash-Williams.
\newblock Edge-disjoint spanning trees of finite graphs.
\newblock {\em J.\ London Math.\ Soc.}, 36:445--450, 1961.

\bibitem{NW71}
C.St.J.A. Nash-Williams.
\newblock Hamiltonian lines in infinite graphs with few vertices of small
  valency.
\newblock {\em Aequationes {M}ath.}, 7:59--81, 1971.

\bibitem{oberlySumner}
D.J. Oberly and D.P. Sumner.
\newblock Every connected, locally connected nontrivial graph with no induced
  claw is hamiltonian.
\newblock {\em J.~Graph Theory}, 3:351--356, 1979.

\bibitem{Oxley79}
J.G. Oxley.
\newblock On a packing problem for infinite graphs and independence.
\newblock {\em J.~Combin.\ Theory (Series B)}, 26:123--130, 1979.

\bibitem{minim}
M.~Stein.
\newblock Ends and vertices of small degree in infinite minimally
  $k$-(edge)-connected graphs.
\newblock SIAM Journal on Discrete Mathematics, in press.

\bibitem{arbo}
M.~Stein.
\newblock Arboricity and tree-packing in locally finite graphs.
\newblock {\em J.~Combin.\ Theory (Series B)}, 96:302--312, 2006.

\bibitem{hcs}
M.~Stein.
\newblock Forcing highly connected subgraphs.
\newblock {\em J.~Graph Theory}, 54:331--349, 2007.

\bibitem{mayaMCW2010}
M.~Stein.
\newblock Degree and substructure in infinite graphs.
\newblock KAM series preprint (MCW 2010), 2010.

\bibitem{LCM}
M.~Stein and J.~Zamora.
\newblock Forcing large complete minors in infinite graphs.
\newblock Preprint 2010.

\bibitem{thomason}
A.G. Thomason.
\newblock The extremal function for complete minors.
\newblock {\em J. Combin. Theory B}, 81:318--338, 2001.

\bibitem{thomassenNonSep}
C.~Thomassen.
\newblock Nonseparating cycles in $k$-connected graphs.
\newblock {\em J.~Graph Theory}, 5:351--354, 1981.

\bibitem{ThomGraphDec}
C.~Thomassen.
\newblock Graph decompositions with applications to subdivisions and path
  systems modulo $k$.
\newblock {\em J. Graph Theory}, 7:261--271, 1983.

\bibitem{tutte61}
W.T. Tutte.
\newblock On the problem of decomposing a graph into {$n$} connected factors.
\newblock {\em J.\ London Math.\ Soc.}, 36:221--230, 1961.

\bibitem{yuV}
X.~Yu.
\newblock Infinite paths in planar graphs {V}, 3-indivisible graphs.
\newblock {\em J. Graph Theory}, 57(4):275--312, 2008.

\end{thebibliography}

\small
\vskip1.8cm 

\noindent Maya Stein\\
{\tt <mstein@dim.uchile.cl>}\\
Centro de Modelamiento Matem\' atico\\
Universidad de Chile\\ 
Blanco Encalada, 2120\\ 
Santiago\\
Chile

\end{document}